\documentclass[dvipdfmx]{amsart}
\usepackage{amsmath}
\usepackage{amssymb} 
\usepackage{amscd} 
\usepackage{graphicx} 
\usepackage{epsfig}
\usepackage[dvips]{color}

\newtheorem{theorem}{Theorem}[section]
\newtheorem{lemma}[theorem]{Lemma}

\newtheorem{proposition}[theorem]{Proposition}
\newtheorem{corollary}[theorem]{Corollary}
\newtheorem{claim}[theorem]{Claim}

\newtheorem{remark}[theorem]{Remark}
\newtheorem{question}[theorem]{Question}
\newtheorem{problem}[theorem]{Problem}

\theoremstyle{definition}
\newtheorem{definition}[theorem]{Definition}
\newtheorem{example}[theorem]{Example}

\numberwithin{equation}{section}
\numberwithin{figure}{section}
\numberwithin{table}{section}

\begin{document}
\baselineskip 14pt

\title{L-space surgery and twisting operation}

\author[K. Motegi]{Kimihiko Motegi}
\address{Department of Mathematics, Nihon University, 
3-25-40 Sakurajosui, Setagaya-ku, 
Tokyo 156--8550, Japan}
\email{motegi@math.chs.nihon-u.ac.jp}

\dedicatory{}

\begin{abstract}
A knot in the $3$--sphere is called an L-space knot if it admits a nontrivial Dehn surgery yielding an L-space, 
i.e. a rational homology $3$--sphere with the smallest possible Heegaard Floer homology.  
Given a knot $K$, 
take an unknotted circle $c$ and twist $K$ $n$ times along $c$ to obtain a twist family $\{ K_n \}$. 
We give a sufficient condition for $\{ K_n \}$ to contain infinitely many L-space knots. 
As an application we show that for each torus knot and each hyperbolic Berge knot $K$, 
we can take $c$ so that the twist family $\{ K_n \}$ contains infinitely many hyperbolic L-space knots. 
We also demonstrate that there is a twist family of hyperbolic L-space knots each member of which has tunnel number greater than one. 
\end{abstract}

\maketitle

{
\renewcommand{\thefootnote}{}
\footnotetext{ 
}

\section{Introduction}
\label{introduction}

Heegaard Floer theory (with $\mathbb{Z}/2\mathbb{Z}$ coefficients) associates  
a group $\widehat{\mathrm{HF}}(M, \textbf{t})$ to a closed, orientable $\mathrm{spin}^c$ 
$3$--manifold $(M, \mathbf{t})$. 
The direct sum of $\widehat{\mathrm{HF}}(M, \textbf{t})$ for all 
$\mathrm{spin}^c$ structures is denoted by $\widehat{\mathrm{HF}}(M)$. 
A rational homology $3$--sphere $M$ is called an \textit{L-space} if 
$\widehat{\mathrm{HF}}(M, \textbf{t})$ is isomorphic to $\mathbb{Z}/2\mathbb{Z}$ for all $\mathrm{spin}^c$ structure $\textbf{t} \in Spin^c(M)$. 
Equivalently, the dimension $\mathrm{dim}_{\mathbb{Z}/2\mathbb{Z}}\widehat{\mathrm{HF}}(M)$ is equal to the order $|H_1(M; \mathbb{Z})|$. 
A knot $K$ in the $3$--sphere $S^3$ is called an  \textit{L-space knot} if the result $K(r)$ of $r$--surgery on $K$ is an L-space for some non-zero integer $r$, 
and the pair $(K, r)$ is called an \textit{L-space surgery}. 
The class of L-spaces includes lens spaces (except $S^2 \times S^1$), 
and more generally, $3$--manifolds with elliptic geometry \cite[Proposition~2.3]{OS3}. 
Since the trivial knot, nontrivial torus knots and Berge knots \cite{Berge2} admit nontrivial surgeries yielding lens spaces, 
these are fundamental examples of L-space knots. 
For the mirror image $K^*$ of $K$, 
$K^*(-r)$ is orientation reversingly homeomorphic to $K(r)$. 
So if $K(r)$ is an L-space, then $K^*(-r)$ is also an L-space \cite[p.1288]{OS3}. 
Hence if $K$ is an L-space knot, 
then so is $K^*$. 

Let $K$ be a nontrivial L-space knot with a positive L-space surgery, 
then Ozsv\'ath and Szab\'o  \cite[Proposition 9.6]{OS4}  (\cite[Lemma 2.13]{Hedden}) prove that 
$r$--surgery on $K$ results in an L-space if and only if $r \ge 2g(K) - 1$, 
where $g(K)$ denotes the genus of $K$. 
This result, 
together with Thurston's hyperbolic Dehn surgery theorem \cite{T1, T2, BePe, PetPorti, BoileauPorti}, 
shows that each hyperbolic L-space knot, 
say a hyperbolic Berge knot, 
produces infinitely many hyperbolic L-spaces by Dehn surgery. 

On the other hand, 
there are some strong constraints for L-space knots: 

\smallskip

$\bullet$ The non-zero coefficients of the Alexander polynomial of an L-space knot are $\pm 1$ and alternate in sign \cite[Corollary~1.3]{OS3}. 

$\bullet$ An L-space knot is fibered \cite[Corollary~1.2]{Ni}(\cite{Ni2}); see also \cite{Ghig, Juh}. 

$\bullet$ An L-space knot is prime \cite[Theorem~1.2]{Kr}. 

\medskip

Note that these conditions are not sufficient. 
For instance, $10_{132}$ satisfies the above conditions, 
but it is not an L-space knot; see \cite{OS3}. 

\medskip

As shown in \cite{Hedden, HLV}, 
some satellite operations keep the property of being L-space knots. 
In the present article, 
we consider if some suitably chosen twistings also keep the property of being L-space knots. 
Given a knot $K$, take an unknotted circle $c$ which bounds a disk intersecting $K$ at least twice.  
Then performing $n$--twist, 
i.e. $(-1/n)$--surgery along $c$, 
we obtain another knot $K_n$. 
Then our question is formulated as: 

\begin{question}
\label{motivation}
Which knots $K$ admit an unknotted circle $c$ such that $n$--twist along $c$ converts $K$ into an L-space knot $K_n$ for infinitely many integers $n$?  
Furthermore, if $K$ has such a circle $c$, 
which circles enjoy the desired property?
\end{question} 

\begin{example}
\label{Pretzel}
Let $K$ be a pretzel knot $P(-2, 3, 7)$ and take an unknotted circle $c$ as in Figure~\ref{pretzel}. 
Then following Ozsv\'ath and Szab\'o \cite{OS3} $K_n$ is an L-space knot if $n \ge -3$ and 
thus the twist family $\{ K_n \}$ contains infinity many L-space knots. 
Note that this family, together with a twist family $\{ T_{2n+1, 2} \}$, comprise all Montesinos L-space knots; 
see \cite{LM} and \cite{BM}. 

\end{example}

\begin{figure}[htbp]
\begin{center}
\includegraphics[width=0.25\linewidth]{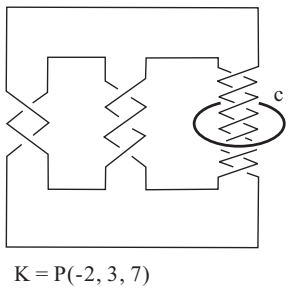}
\caption{A knot $K_n$ obtained by $n$--twist along $c$ is an L-space knot if $n \ge -3$.}
\label{pretzel}
\end{center}
\end{figure}

In this example, 
it turns out that $c$ becomes a Seifert fiber in the lens space $K(19)$
(cf. Example~\ref{T53_2_nontorus}).  
We employed such a circle for relating Seifert fibered surgeries in \cite{DMM1}. 
A pair $(K, m)$ of a knot $K$ in $S^3$ and an integer $m$ is 
a \textit{Seifert surgery} if $K(m)$ has a Seifert fibration;
we allow the fibration to be degenerate,
i.e.\ it contains an exceptional fiber of index 0 as a degenerate fiber. 
See \cite[2.1]{DMM1} for details. 
The definition below enables us to say that $c$ is a seiferter for the Seifert (lens space) surgery $(K, 19)$. 

\begin{definition}[\textbf{seiferter} \cite{DMM1}]
\label{seiferter}
Let $(K, m)$ be a Seifert surgery.
A knot $c$ in $S^3 -N(K)$ is called a \textit{seiferter} for $(K, m)$
if $c$ satisfies the following: 
\begin{itemize}
\item $c$ is a trivial knot in $S^3$.
\item $c$ becomes a fiber in
a Seifert fibration of $K(m)$.
\end{itemize}
\end{definition}

As remarked in \cite[Convention~2.15]{DMM1}, 
if $c$ bounds a disk in $S^3 - K$, 
then we do not regard $c$ as a seiferter. 
Thus for any seiferter $c$ for $(K, m)$, 
$S^3 - \mathrm{int}N(K \cup c)$ is irreducible.

Let $(K, m)$ be a Seifert surgery with a seiferter $c$. 
There are two cases according as $c$ becomes a fiber in a non-degenerate Seifert fibration 
of $K(m)$ or $c$ becomes a fiber in a degenerate Seifert fibration of $K(m)$. 
In the former case, 
for homological reasons, the base surface is 
the $2$--sphere $S^2$ or the projective plane 
$\mathbb{R}P^2$. 
Suppose that $c$ is a fiber in a non-degenerate Seifert fibration of $K(m)$ over the $2$--sphere $S^2$. 
Then in the following we assume that $K(m)$ contains at most three exceptional fibers and 
if there are three exceptional fibers, 
then $c$ is an exceptional fiber. 
We call such a seiferter a \textit{seiferter for a small Seifert fibered surgery $(K, m)$}. 
To be precise, 
the images of $K$ and $m$ after $n$--twist along $c$ should be denoted by $K_{c, n}$ and $m_{c, n}$, 
but for simplicity, 
we abbreviate them to $K_n$ and $m_n$ respectively as long as there is no confusion. 

\begin{theorem}
\label{tree small}
Let $c$ be a seiferter for a small Seifert fibered surgery $(K, m)$. 
Then $(K_n, m_n)$ is an L-space surgery for an infinite interval of integers $n$ 
if and only if the result of $(m, 0)$--surgery on $K \cup c$ is an L-space.   
\end{theorem}

In remaining cases, 
it turns out that every seiferter enjoys the desired property in Question~\ref{motivation}. 

\begin{theorem}
\label{tree projective}
Let $c$ be a seiferter for $(K, m)$ which become a fiber in a Seifert fibration of $K(m)$ over $\mathbb{R}P^2$.
Then $(K_n, m_n)$ is an L-space surgery for all but at most one integer $n_0$ with 
$(K_{n_0}, m_{n_0}) = (O, 0)$. 
Hence $K_n$ is an L-space knot for all integers $n$. 
\end{theorem}

Let us turn to the case where $c$ is a $($degenerate or non-degenerate$)$ fiber in a degenerate Seifert fibration of $K(m)$.  
Recall from \cite[Proposition~2.8]{DMM1} that 
if $K(m)$ has a degenerate Seifert fibration, 
then it is a lens space or a connected sum of two lens spaces such that each summand is neither $S^3$ nor $S^2 \times S^1$. 
The latter $3$--manifold will be simply referred to as a \textit{connected sum of two lens spaces}, 
which is an L-space \cite[8.1(5)]{Szabo} (\cite{OS2}). 

\begin{theorem}
\label{tree degenerate}
Let $c$ be a seiferter for $(K, m)$ which becomes a $($degenerate or non-degenerate$)$ fiber in 
a degenerate Seifert fibration of $K(m)$. 
\begin{enumerate}
\item
If $K(m)$ is a lens space, 
then $(K_n, m_n)$ is an L-space surgery, hence $K_n$ is an L-space knot, 
for all but at most one integer $n$. 
\item
If $K(m)$ is a connected sum of two lens spaces, 
then $(K_n, m_n)$ is an L-space surgery, 
hence $K_n$ is an L-space knot, 
for any integer $n \ge -1$ or $n \le 1$. 
\end{enumerate}
\end{theorem}

Following Greene \cite[Theorem~1.5]{Greene}, 
if $K(m)$ is a connected sum of two lens spaces, 
then $K$ is a torus knot $T_{p, q}$ or a cable of a torus knot $C_{p, q}(T_{r, s})$, where $p = qrs \pm 1$.  
We may assume $p, q \ge 2$ by taking the mirror image if necessary. 
The next theorem is a refinement of Theorem~\ref{tree degenerate}(2). 

\begin{theorem}
\label{Tpq}
Let $c$ be a seiferter for $(K, m) = (T_{p, q}, pq)$ or $(C_{p, q}(T_{r, s}), pq)$ 
$(p = qrs \pm 1)$. 
We assume $p, q \ge 2$. 
Then a knot $K_n$ obtained from $K$ by $n$--twist along $c$ 
is an L-space knot for any $n \ge -1$. 
Furthermore, 
if the linking number $l$ between $c$ and $K$ satisfies $l^2 \ge 2pq$, 
then $K_n$ is an L-space knot for all integers $n$. 
\end{theorem}

In the above theorem, 
even when $\ell^2 < 2pq$, 
$K_n$ ($n < -1$) may be an L-space knot; see \cite{MT}. 

In Sections~\ref{twisted torus knots}, \ref{twisted Berge knots} and \ref{twisted unknots} 
we will exploit seiferter technology developed in \cite{DMM1, DMM2, DEMM} to give 
a partial answer to Question~\ref{motivation}. 
Even though Theorem~\ref{Tpq} treats a special kind of Seifert surgeries, 
it offers many applications. 
In particular, Theorem~\ref{Tpq} enables us to give new families of L-space twisted torus knots. 
See Section~\ref{twisted torus knots} for the definition of twisted torus knots $K(p, q; r, n)$  introduced by Dean \cite{D}. 

\begin{theorem}[\textbf{L-space twisted torus knots}]
\label{twisted torus knot}
\begin{enumerate}
\item
The following twisted torus knots are L-space knots for all integers $n$. 
	\begin{itemize}
	\item
	$K(p, q; p + q, n)$ with $p, q \ge 2$
	\item
	$K(3p+1, 2p+1; 4p+1, n)$ with $p > 0$
	\item
	$K(3p+2, 2p+1; 4p+3, n)$ with $p > 0$ 
	\end{itemize}

\item
The following twisted torus knots are L-space knots for any $n \ge -1$. 
	\begin{itemize}
	\item
	$K(p, q; p - q, n)$ with $p, q \ge 2$ 
	\item
	$K(2p+3, 2p+1; 2p+2, n)$ with $p > 0$
	\end{itemize}
\end{enumerate}
\end{theorem}

Theorem~\ref{twisted torus knot} has the following corollary, 
which asserts that every nontrivial torus knot admits twistings desired in Question~\ref{motivation}. 

\begin{corollary}
\label{L-space twisted torus knot} 
For any nontrivial torus knot $T_{p, q}$, 
we can take an unknotted circle $c$ so that $n$--twist along $c$ 
converts $T_{p, q}$ into an L-space knot $K_n$ for all integers $n$. 
Furthermore, $\{ K_n\}_{|n| > 3}$ is a set of mutually distinct hyperbolic L-space knots. 
\end{corollary}

For the simplest L-space knot, 
i.e. the trivial knot $O$,  
we can strengthen Corollary~\ref{L-space twisted torus knot} as follows. 

\begin{theorem}[\textbf{L-space twisted unknots}]
\label{L-space twisted unknot}
For the trivial knot $O$, 
we can take infinitely many unknotted circles $c$ so that 
$n$--twist along $c$ changes $O$ into a nontrivial L-space knot $K_{c, n}$ for any non-zero integer $n$. 
Furthermore, $\{ K_{c, n}\}_{|n| > 1}$ is a set of mutually distinct hyperbolic L-space knots. 
\end{theorem}

Using a relationship between Berge's lens space surgeries and surgeries yielding a connected sum of two lens spaces,  
we can prove: 

\begin{theorem}[\textbf{L-space twisted Berge knots}]
\label{Berge}
For any hyperbolic Berge knot $K$, 
there is an unknotted circle $c$ such that 
$n$--twist along $c$ converts $K$ into a hyperbolic L-space knot $K_n$ for infinitely many 
integers $n$. 
\end{theorem}

In Section~\ref{tunnel} we consider the tunnel number of L-space knots. 
Recall that the \textit{tunnel number} of a knot $K$ in $S^3$ is the minimum number of mutually disjoint, 
embedded arcs connecting $K$ such that the exterior of the resulting $1$--complex is a handlebody. 
Hedden's cabling construction \cite{Hedden}, 
together with \cite{MoSa},  
enables us to obtain an L-space knot with tunnel number greater than $1$. 
Actually Baker and Moore \cite{BM} have shown that for any integer $N$, 
there is an L-space knot with tunnel number greater than $N$. 
However,  L-space knots with tunnel number greater than one constructed above are all satellite (non-hyperbolic) knots and they ask: 

\begin{question}[\cite{BM}]
\label{BMquestion}
Is there a non-satellite, L-space knot with tunnel number greater than one? 
\end{question}

Examining knots with Seifert surgeries which do not arise from primitive/Seifert-fibered construction given by \cite{EJMM}, 
we prove the following which answers the question in the positive. 

\begin{theorem}
\label{tunnel2}
There exist infinitely many hyperbolic L-space knots with tunnel number greater than one.   
\end{theorem}

Each knot in the theorem is obtained from a trefoil knot $T_{3, 2}$ by alternate twisting along two seiferters 
for the lens space surgery $(T_{3, 2}, 7)$. 

In Section~\ref{questions} we will discuss further questions on relationships between L-space knots and twisting operation. 

\bigskip

\noindent
\textbf{Acknowledgments.}
I would like to thank Ken Baker for insightful and stimulating discussion.  
This paper would not have come to fruition without examples obtained in previous joint works 
with Arnaud Deruelle, Mario Eudave-Mu\~noz, Edgar Jasso and Katura Miyazaki. 
Finally, I would like to thank the referee for careful reading and useful comments. 

The author has been partially supported by JSPS Grants--in--Aid for Scientific 
Research (C), 26400099, The Ministry of Education, Culture, Sports, Science and Technology, Japan and Joint Research Grant of Institute of Natural Sciences at Nihon University for 2014. 

\bigskip

\section{Seifert fibered L-spaces}
\label{SeifertL}
Let $M$ be a rational homology $3$--sphere which is a Seifert fiber space. 
For homological reasons, the base surface of $M$ is either $S^2$ or $\mathbb{R}P^2$. 
In the latter case, 
Boyer, Gordon and Watson \cite[Proposition~5]{BGW} prove that $M$ is an L-space. 
Now assume that the base surface of $M$ is $S^2$.  
Following Ozsv\'ath and Szab\'o \cite[Theorem~1.4]{OSgenus} 
if $M$ is an L-space, 
then it carries no taut foliation, 
in particular, it carries no horizontal (i.e. transverse) foliation. 
Furthermore, 
Lisca and Stipsicz \cite[Theorem~1.1]{LS} prove that the converse also does hold.  
Therefore a Seifert fibered rational homology $3$--sphere $M$ over $S^2$ is an L-space if and only if 
it does not admit a horizontal foliation.  
Note that if $M$ does not carry a horizontal foliation, 
then it is necessarily a rational homology $3$--sphere. 
In fact, if $|H_1(M; \mathbb{Z})| = \infty$, 
then $M$ is a surface bundle over the circle \cite[VI.34]{Ja}, \cite{Hat}, 
and hence it has a horizontal foliation. 
On the other hand, Eisenbud-Hirsh-Neumann \cite{EHN}, 
Jankins-Neumann \cite{JN} and Naimi \cite{Na} gave a necessary and sufficient
conditions for a Seifert fibered 3-manifold to carry a horizontal foliation. 
Combining them 
we have Theorem~\ref{SeifertLspace} below. 
See also \cite[Theorem 5.4]{BRW}; we follow the convention of Seifert invariants in \cite[Section 4]{BRW}. 

For ordered triples $(a_1, a_2, a_3)$ and $(b_1, b_2, b_3)$, 
we write $(a_1, a_2, a_3) < (b_1, b_2, b_3)$ (resp. $(a_1, a_2, a_3) \le (b_1, b_2, b_3)$) 
if $a_i < b_i$ (resp. $a_i \le b_i$) for $1 \le i \le 3$, 
and denote by $(a_1, a_2, a_3)^*$ the ordered triple 
$(\sigma(a_1), \sigma(a_2), \sigma(a_3))$, 
where $\sigma$ is a permutation such that 
$\sigma(a_1) \le \sigma(a_2) \le \sigma(a_3)$. 

\begin{theorem}[\cite{OSgenus, LS, EHN, JN, Na}]
\label{SeifertLspace}
A Seifert fiber space $S^2(b, r_1, r_2, r_3)$ $(b \in \mathbb{Z},\ 0 < r_i < 1)$ is 
an L-space if and only if one of the following holds. 
\begin{enumerate}
\item
$b  \ge 0$ or $b \le -3$. 
\item
$b = -1$ and there are no relatively prime integers $a, k$ such that $0 < a \le k/2$ and 
$(r_1, r_2, r_3)^* < (1/k, a/k, (k-a)/k)$. 
\item
$b = -2$ and there are no relatively prime integers $0 < a \le k/2$ such that 
$(1-r_1, 1-r_2, 1-r_3)^* < (1/k, a/k, (k-a)/k)$. 
\end{enumerate}
\end{theorem}

For our purpose, we consider the following problem: 

\begin{problem}
\label{solution}
Given an integer $b$ and rational numbers $0 < r_1 \le r_2 < 1$, 
describe rational numbers $-1 \le r \le 1$ for which $S^2(b, r_1, r_2, r)$ is an L-space. 
\end{problem}

We begin by observing: 

\begin{lemma}
\label{any r}
Assume that $0 < r_1 \le r_2 < 1$. 
\begin{enumerate}
\item
If $b \ge 0$ or $b \le -3$, 
then $S^2(b, r_1, r_2, r)$ is an L-space for any $0 < r < 1$. 
\item
If $r_1 + r_2 \ge 1$, 
then $S^2(-1, r_1, r_2, r)$ is an L-space for any $0 < r < 1$. 
\item
If $r_1 + r_2 \le 1$, 
then  $S^2(-2, r_1, r_2, r)$ is an L-space for any $0 < r < 1$. 
\end{enumerate}
\end{lemma}

\noindent
\textit{Proof of Lemma~\ref{any r}.}
The first assertion is nothing but Theorem~\ref{SeifertLspace}(1). 

Suppose for a contradiction that 
$S^2(-1, r_1, r_2, r)$ is not an L-space for some $0 < r < 1$.
Then, by Theorem~\ref{SeifertLspace}(2) we can take relatively prime integers $a, k$
$(0 < a \le k/2)$ so that 
$(r_1, r_2, r)^* < (1/k, a/k, (k-a)/k)$. 
This then implies that 
$r_ 1 < a/k$ and $r_2 < (k-a)/k$. 
Hence $r_1 + r_2 < a/k + (k-a)/k = 1$, a contradiction. 
This proves $(2)$.  

To prove $(3)$, assume for a contradiction that 
$S^2(-2, r_1, r_2, r)$ is not an L-space for some $0 < r < 1$. 
Then, 
by Theorem~\ref{SeifertLspace}(3) we have relatively prime integers $a, k$ 
$(0 < a \le k/2)$ such that 
$(1-r_1, 1-r_2, 1-r)^* < (1/k, a/k, (k-a)/k)$. 
Thus we have 
$(1-r_2) < a/k$ and  $(1-r_1) < (k-a)/k$. 
Thus $(1-r_1) + (1-r_2) < 1$, 
which implies $r_1 + r_2 > 1$, 
contradicting the assumption. 
\hspace*{\fill} $\square$(Lemma~\ref{any r})

\medskip

Now let us prove the following, which gives an answer to Problem~\ref{solution}. 

\begin{proposition}
\label{r_c=0}
Assume that $0 < r_1 \le r_2 < 1$. 
\begin{enumerate}
\item
If $b \le -3$ or $b \ge 1$, 
then $S^2(b, r_1, r_2, r)$ is an L-space for any $-1 \le r \le 1$.  

\item
If $b = -2$, 
then there exists $\varepsilon > 0$ such that 
$S^2(-2, r_1, r_2, r)$ is an L-space for any $-1 \le r \le \varepsilon$. 
Furthermore, if $r_1 + r_2 \le 1$, 
then $S^2(-2, r_1, r_2, r)$ is an L-space if $-1 \le r < 1$. 

\item
Suppose that $b = -1$. 
	\begin{enumerate}
	\item
	If $r _1 + r_2 \ge 1$, 
	then $S^2(-1, r_1, r_2, r)$ is an L-space for any $0 < r \le 1$.  
	\item 
	If $r _1 + r_2 \le 1$, 
	then $S^2(-1, r_1, r_2, r)$ is an L-space for any $-1 \le r  < 0$. 
	\end{enumerate}

\item
If $b = 0$, 
then there exists $\varepsilon > 0$ such that 
$S^2(r_1, r_2, r)$ is an L-space for any $-\varepsilon \le r \le 1$. 
Furthermore, if $r_1 + r_2 \ge 1$, 
then $S^2(r_1, r_2, r)$ is an L-space if $-1 < r \le 1$. 
\end{enumerate}
\end{proposition}

\noindent
\textit{Proof of Proposition~\ref{r_c=0}.}
If $r= 0, \pm 1$, then $S^2(b, r_1, r_2, r)$ is a lens space. 

\begin{claim}
\label{S2xS1}
Suppose that $r$ is an integer. 
Then the lens space $S^2(b, r_1, r_2, r)$ is $S^2 \times S^1$ if and only if 
$b+r =-1$ and $r_1+r_2 = 1$. 
In particular, 
if $b+r \ne -1$, then $S^2(b, r_1, r_2, r)$ is an L-space. 
\end{claim}

\noindent
\textit{Proof of Claim~\ref{S2xS1}.}
Recall that $H_1(S^2(a/b, c/d)) \cong \mathbb{Z}$ $(b, d \ge 1)$ if and only if 
$ad + bc = 0$, i.e. $a/b + c/d = 0$. 
Thus $S^2(b, r_1, r_2, r)$ is $S^2 \times S^1$ if and only if 
$b + r_1 + r_2 + r =0$, i.e. $r_1 + r_2 = -b-r \in \mathbb{Z}$. 
Since $0 < r_i < 1$, we have $r_1+ r_2 = 1$ and $b+r = -1$. 
\hspace*{\fill} $\square$(Claim~\ref{S2xS1})

\medskip

We divide into two cases according as $0 \le r \le 1$ or $-1 \le  r \le  0$. 

\medskip

\noindent
\textbf{Case I}.\quad $0 \le r \le 1$. 

\smallskip

(i) If $b \ge 0$ or $b \le -3$, 
then $S^2(b, r_1, r_2, r)$ is an L-space for any $0 < r < 1$ by Lemma~\ref{any r}(1).  
Since $b+r \ne -1$ for $r  = 0, 1$, 
by Claim~\ref{S2xS1} $S^2(b, r_1, r_2, r)$ is an L-space for $r = 0, 1$. 
Hence $S^2(b, r_1, r_2, r)$ is an L-space for any $0 \le r \le 1$. 

\smallskip 

(ii) Suppose that $b=-1$.  
By Lemma~\ref{any r}(2), if $r_1 + r_2 \ge 1$, 
then $S^2(-1, r_1, r_2, r)$ is an L-space for any $0 < r < 1$. 
Since $S^2(-1, r_1, r_2, 1)$ is an L-space (Claim~\ref{S2xS1}), 
$S^2(-1, r_1, r_2, r)$ is an L-space for any $0 < r \le 1$. 

\smallskip 

(iii) Assume $b= -2$. 
Let us assume $0 < r \le r_1$ so that 
$0 < 1-r_2 \le 1-r_1 \le 1-r < 1$. 
Set $A = \{ (k-a)/k\ |\ 1-r_2 < 1/k,\ 1-r_1 < a/k,\ 0 < a \le k/2,\ \textrm{$a$ and $k$ are relatively prime integers} \}$. 
If $A = \emptyset$, 
i.e. there are no relatively prime integers $a, k$ $(0 < a \le k/2)$ such that $1-r_2 < 1/k,\ 1-r_1 < a/k$, 
then $S^2(-2, r_1, r_2, r)$ is an L-space for any $0 < r \le r_1$ by Theorem~\ref{SeifertLspace}. 
Suppose that $A \ne \emptyset$. 
Since there are only finitely many integers $k$ satisfying $1-r_2 < 1/k$, 
$A$ consists of only finitely many elements. 
Let $r_0$ be the maximal element in $A$. 
If $0 < r \le 1-r_0$, then $r_0 \le 1-r <1$, 
and hence there are no relatively prime integers $a, k$ $(0 < a \le k/2)$ satisfying 
$(1-r_2, 1-r_1, 1-r) < (1/k, a/k, (k-a)/k)$. 
Put $\varepsilon = \mathrm{min}\{ r_1, 1-r_0 \}$. 
Then $S^2(-2, r_1, r_2, r)$ is an L-space for any $0 < r \le \varepsilon$ by Theorem~\ref{SeifertLspace}. 
Since $S^2(-2, r_1, r_2, 0)$ is an L-space (Claim~\ref{S2xS1}), 
$S^2(-2, r_1, r_2, r)$ is an L-space for any $0 \le r \le \varepsilon$. 
Furthermore, 
if we have the additional condition $r_1 + r_2 \le 1$, 
then Lemma~\ref{any r}(3) improves the result so that 
$S^2(-2, r_1, r_2, r)$ is an L-space for any $0 \le r < 1$.

\medskip

\noindent
\textbf{Case II}. \quad $-1 \le r \le 0$. 

Note that $S^2(b, r_1, r_2, r) = S^2(b-1, r_1, r_2, r+1)$. 

\smallskip

(i) If $b \ge 1$ or $b \le -2$ (i.e. $b-1 \ge 0$ or $b-1 \le -3$), 
then $S^2(b, r_1, r_2, r) = S^2(b-1, r_1, r_2, r+1)$ is an L-space for any 
$0 < r+1 < 1$, i.e. $-1 < r < 0$ by Lemma~\ref{any r}(1). 
Since $b+r \ne -1$ for $r = -1, 0$,  
$S^2(b, r_1, r_2, r)$ is an L-space for $r = -1, 0$ (Claim~\ref{S2xS1}). 
Thus $S^2(b, r_1, r_2, r)$ is an L-space for any $-1 \le r \le 0$. 

\smallskip

(ii) If $b=0$ (i.e. $b-1 = -1$), 
then $S^2(0, r_1, r_2, r) = S^2(-1, r_1, r_2, r+1)$. 
Let us assume $r_2 -1 \le r < 0$ so that 
$0 < r_1 \le r_2 \le r+1 < 1$. 
Set $A = \{ (k-a)/k \ |\ r_1< 1/k,\ r_2 < a/k,\ 0 < a \le k/2,\ \textrm{$a$ and $k$ are relatively prime integers} \}$. 
If $A = \emptyset$, 
then we can easily observe that for any $r$ with $r_2 \le r+1 < 1$, 
$S^2(-1, r_1, r_2, r+1)$ is an L-space (Theorem~\ref{SeifertLspace}). 
Hence for any $r_2 -1 \le r < 0$, 
$S^2(0, r_1, r_2, r)$ is an L-space. 
Suppose that $A \ne \emptyset$. 
Since $A$ is a finite set, 
we take the maximal element $r_0$ in $A$. 
If $r_0 \le r+1 < 1$ (i.e. $r_0-1 \le r < 0$), 
then there are no relatively prime integers $a, k$ $(0 < a \le k/2)$ satisfying 
$(r_1, r_2, r+1) < (1/k, a/k, (k-a)/k)$. 
Put $\varepsilon = \mathrm{min}\{ 1-r_2,\ 1-r_0 \}$. 
Then $S^2(0, r_1, r_2, r) = S^2(-1, r_1, r_2, r+1)$ is an L-space for any $-\varepsilon \le r < 0$ (Theorem~\ref{SeifertLspace}). 
Since $S^2(0, r_1, r_2, 0) = S^2(r_1, r_2)$ is an L-space (Claim~\ref{S2xS1}), 
$S^2(0, r_1, r_2, r)$ is an L-space for any $-\varepsilon \le r \le 0$.   
Furthermore, 
if we have the additional condition $r_1 + r_2 \ge 1$, 
then Lemma~\ref{any r}(2) improves the result so that 
$S^2(r_1, r_2, r) = S^2(-1, r_1, r_2, r+1)$ is an L-space for any $-1 < r \le 0$. 

\smallskip

(iii) If $b = -1$ (i.e. $b-1 = -2$), 
then $S^2(-1, r_1, r_2, r) = S^2(-2, r_1, r_2, r+1)$.  
Assume that $r_1 + r_2 \le 1$. 
Then Proposition~\ref{any r}(3), 
$S^2(-1, r_1, r_2, r) = S^2(-2, r_1, r_2, r+1)$ is an L-space for any 
$0 < r+1 < 1$, i.e. $-1 < r < 0$. 
Since Claim~\ref{S2xS1} shows that $S^2(-1, r_1, r_2, -1)$ is an L-space, 
$S^2(-1, r_1, r_2, r)$ is an L-space for any 
$-1 \le r < 0$. 

\medskip
 
Combining Cases I and II, we obtain the result described in the proposition. 

\hspace*{\fill} $\square$(Proposition~\ref{r_c=0})

\medskip

The next proposition shows that if $S^2(b, r_1, r_2, r_{\infty})$ is an L-space for some rational number 
$0 < r_{\infty} < 1$, 
then we can find $r$ near $r_{\infty}$ so that $S^2(b, r_1, r_2, r)$ is an L-space. 

\begin{proposition}
\label{b=-1, -2}
Suppose that $0 < r_1 \le r_2 < 1$ and 
$S^2(b, r_1, r_2, r_{\infty})$ is an L-space for some rational number 
$0 < r_{\infty} < 1$. 
\begin{enumerate}
\item
If $b = -1$, 
then $S^2(-1, r_1, r_2, r)$ is an L-space for any $r_{\infty} \le r \le 1$. 
\item
If $b = -2$, 
then $S^2(-2, r_1, r_2, r)$ is an L-space for any $-1 \le r \le r_{\infty}$. 
\end{enumerate}
\end{proposition}

\noindent
\textit{Proof of Proposition~\ref{b=-1, -2}.}
$(1)$ Assume for a contradiction that 
$S^2(-1, r_1, r_2, r)$ is not an L-space for some $r$ satisfying $r_{\infty} \le r < 1$. 
By Theorem~\ref{SeifertLspace} we have relatively prime integers $a, k$ $(0 < a \le k/2)$ such that 
$(r_1, r_2, r)^* < (1/k, a/k, (k-a)/k)$. 
Since $r_{\infty} \le r < 1$,  
$(r_1, r_2, r_{\infty})^* \le (r_1, r_2, r)^* < (1/k, a/k, (k-a)/k)$. 
Hence Theorem~\ref{SeifertLspace} shows that $S^2(-1, r_1, r_2, r_{\infty})$ is not an L-space, 
a contradiction.  
Since $S^2(-1, r_1, r_2, 1) = S^2(r_1, r_2)$ is an L-space (Claim~\ref{S2xS1}), 
$S^2(-1, r_1, r_2, r)$ is an L-space for any $r_{\infty} \le r \le 1$. 

$(2)$ Next assume for a contradiction that $S^2(-2, r_1, r_2, r)$ is not an L-space for some $r$ satisfying 
$0 < r \le r_{\infty}$. 
Then following Theorem~\ref{SeifertLspace} we have 
$(1-r_1, 1-r_2, 1-r)^* < (1/k, a/k, (k-a)/k)$ for some relatively prime integers $a, k$ $(0 < a \le k/2)$. 
Since $r \le r_{\infty}$, we have $1 - r_{\infty} \le 1-r$, 
and hence $(1-r_1, 1-r_2, 1-r_{\infty})^* \le (1-r_1, 1-r_2, 1-r)^* < (1/k, a/k, (k-a)/k)$. 
This means $S^2(-2, r_1, r_2, r_{\infty})$ is not an L-space, contradicting the assumption. 
Thus $S^2(-2, r_1, r_2, r)$ is an L-space for any $0 < r \le r_{\infty}$. 
Furthermore, as shown in Proposition~\ref{r_c=0}(2), 
$S^2(-2, r_1, r_2, r)$ is an L-space if $-1 \le r \le \varepsilon$ for some $\varepsilon > 0$, 
so $S^2(-2, r_1, r_2, r)$ is an L-space for any $-1 \le r \le r_{\infty}$. 
\hspace*{\fill} $\square$(Proposition~\ref{b=-1, -2})

\medskip

We close this section with the following result which is the complement of Proposition~\ref{b=-1, -2}. 

\begin{proposition}
\label{nonLspace}
Suppose that $0 < r_1 \le r_2 < 1$ and 
$S^2(b, r_1, r_2, r_{\infty})$ is not an L-space for some rational number 
$0 < r_{\infty} < 1$. 
\begin{enumerate}
\item
If $b = -1$, 
then there exists $\varepsilon > 0$ such that 
$S^2(-1, r_1, r_2, r)$ is not an L-space for any $0< r < r_{\infty} + \varepsilon$. 
\item
If $b = -2$, 
then there exists $\varepsilon > 0$ such that 
then $S^2(-2, r_1, r_2, r)$ is an L-space for any $r_{\infty}-\varepsilon < r < 1$. 
\end{enumerate}
\end{proposition}

\noindent
\textit{Proof of Proposition~\ref{nonLspace}.}
$(1)$ Since $S^2(-1, r_1, r_2, r_{\infty})$ is not an L-space, 
Theorem~\ref{SeifertLspace} shows that there are relatively prime integers $a, k$ $(0 < a \le k/2)$ such that 
$(r_1, r_2, r_{\infty})^* < (1/k, a/k, (k-a)/k)$. 
Then clearly there exists $\varepsilon > 0$ such that for any $0 < r < r_{\infty}+ \varepsilon$, 
we have $(r_1, r_2, r)^* < (1/k, a/k, (k-a)/k)$. 
Thus by Theorem~\ref{SeifertLspace} again 
$S^2(-1, r_1, r_2, r)$ is not an L-space for any $0 < r < r_{\infty}+ \varepsilon$. 

$(2)$ Since $S^2(-2, r_1, r_2, r_{\infty})$ is not an L-space, 
by Theorem~\ref{SeifertLspace} we have relatively prime integers $a, k$ $(0 < a \le k/2)$ such that 
$(1-r_1, 1-r_2, 1-r_{\infty})^* < (1/k, a/k, (k-a)/k)$. 
Hence there exists $\varepsilon > 0$ such that 
if $0 < 1-r < 1-r_{\infty} + \varepsilon$, 
i.e. $r_{\infty} - \varepsilon < r < 1$, 
then $(1-r_1, 1-r_2, 1-r)^* < (1/k, a/k, (k-a)/k)$. 
Following Theorem~\ref{SeifertLspace} 
$S^2(-2, r_1, r_2, r)$ is not an L-space for any $r_{\infty} - \varepsilon < r < 1$. 
\hspace*{\fill} $\square$(Proposition~\ref{nonLspace})

\bigskip

\section{L-space surgeries and twisting along seiferters I -- non-degenerate case}
\label{Twisting L-space seiferter I}

The goal in this section is to prove Theorems~\ref{tree small} and \ref{tree projective}.  

Let $c$ be a seiferter for a small Seifert fibered surgery $(K, m)$. 
The $3$--manifold obtained by $(m, 0)$--surgery on $K \cup c$ is 
denoted by $M_c(K, m)$. 

\medskip

\noindent
\textit{Proof of Theorem~\ref{tree small}.}
First we prove the ``if" part of Theorem~\ref{tree small}. 
If $K(m)$ is a lens space and $c$ is a core of the genus one Heegaard splitting, 
then $K_n(m_n)$ is a lens space for any integer $n$. 
Thus $(K_n, m_n)$ is an L-space surgery  for all $n \in \mathbb{Z}$ except when $K_n(m_n) \cong S^2 \times S^1$, 
i.e. $K_n$ is the trivial knot and $m_n = 0$ \cite[Theorem~8.1]{GabaiIII}. 
Since $(K_n, m_n) = (K_{n'}, m_{n'})$ if and only if $n = n'$ \cite[Theorem~5.1]{DMM1}, 
there is at most one integer $n$ such that $(K_n, m_n) = (O, 0)$. 
Henceforth, in the case where $K(m)$ is a lens space, 
we assume that $K(m)$ has a Seifert fibration over $S^2$ 
with two exceptional fibers $t_1$, $t_2$, 
and $c$ becomes a regular fiber in this Seifert fibration. 

Let $E$ be $K(m) - \mathrm{int}N(c)$ with a fibered tubular neighborhood of the union of two exceptional fibers $t_1, t_2$
and one regular fiber $t_0$ removed. 
Then $E$ is a product circle bundle over the four times punctured sphere. 
Take a cross section of $E$ 
such that $K(m)$ is expressed as 
$S^2(b, r_1, r_2, r_3)$, 
where the Seifert invariant of $t_0$ is $b \in \mathbb{Z}$, 
that of $t_i$ is $ 0 < r_i < 1$ $(i = 1, 2)$, 
and that of $c$ is $0 \le r_3 < 1$. 
Without loss of generality, we may assume $r_1 \le r_2$. 
Let $s$ be the boundary curve on $\partial N(c)$ of 
the cross section so that $[s]\cdot[t] = 1$ for a regular fiber $t \subset \partial N(c)$. 
Let $(\mu, \lambda)$ be a preferred meridian-longitude pair of $c \subset S^3$. 
Then $[\mu] = \alpha_3 [s] + \beta_3 [t] \in H_1( \partial N(c) )$ and 
$[\lambda] = -\alpha [s] - \beta [t] \in H_1( \partial N(c) )$ for some integers 
$\alpha_3, \beta_3, \alpha$ and $\beta$ which satisfy $\alpha_3 > 0$ and $\alpha \beta_3 - \beta\alpha_3 = 1$, 
where $r_3 = \beta_3 / \alpha_3$. 
Now let us write $r_c = \beta /\alpha$, 
which is the slope of the preferred longitude $\lambda$ of $c \subset S^3$ with respect to $(s, t)$--basis. 

\begin{claim}
\label{Seifert Mc}
$M_c(K, m)$ is a $($possibly degenerate$)$ Seifert fiber space $S^2(b, r_1, r_2, r_c)$; 
if $r_c = -1/0$, 
then it is a connected sum of two lens spaces. 
\end{claim}

\noindent
\textit{Proof of Claim~\ref{Seifert Mc}.}
$M_c(K, m)$ is regarded as a $3$--manifold obtained from $K(m)$ by performing $\lambda$--surgery along the fiber $c \subset K(m)$. 
Since $[\lambda] = -\alpha [s] - \beta [t]$, 
$M_c(K, m)$ is a (possibly degenerate) Seifert fiber space $S^2(b, r_1, r_2, r_c)$. 
If $\alpha = 0$, i.e. $r_c = -1/0$, 
then $M_c(K, m)$ has a degenerate Seifert fibration and it is a connected sum of two lens spaces. 
\hspace*{\fill} $\square$(Claim~\ref{Seifert Mc})

\medskip

Recall that $(K_n, m_n)$ is a Seifert surgery obtained from $(K, m)$ by 
twisting $n$ times along $c$.
The image of $c$ after the $n$--twist along $c$
is also a seiferter for $(K_n, m_n)$ and denoted 
by $c_n$. 
We study how the Seifert invariant of $K(m)$ behaves under the twisting.
We compute the Seifert invariant of
$c_n$ in $K_n(m_n)$ under the same cross section on $E$. 

Since we have 
\begin{eqnarray*}
\begin{pmatrix}
	[\mu] \\[2pt]
	[\lambda]
\end{pmatrix}
=
\begin{pmatrix}
	\alpha_3 & \beta_3 \\[2pt]
	-\alpha & -\beta
\end{pmatrix}
\begin{pmatrix}
	[s] \\[2pt]
	[t]
\end{pmatrix}, 
\end{eqnarray*} 

it follows that

\begin{eqnarray*}
\begin{pmatrix}
	[s] \\[2pt] [t] 
\end{pmatrix}
=
\begin{pmatrix}	
-\beta & -\beta_3 \\[2pt]
	\alpha & \alpha_3 
\end{pmatrix}
        \begin{pmatrix}
	[\mu] \\[2pt]
	[ \lambda ]
\end{pmatrix}.
\end{eqnarray*}

Twisting $n$ times along $c$ is equivalent to
performing $-1/n$--surgery on $c$.
A preferred meridian-longitude pair
$(\mu_n, \lambda_n)$ of $N(c_n) \subset S^3$ 
satisfies
$[\mu_n] = [\mu] - n[\lambda]$ and
$[\lambda_n] = [\lambda]$
in $H_1( \partial N(c_n) ) = H_1(\partial N(c))$.

We thus have 
\begin{eqnarray*}
\begin{pmatrix}
	[s] \\[2pt]
	[t]
\end{pmatrix}
=
\begin{pmatrix}
	-\beta & -n\beta - \beta_3 \\[2pt]
	\alpha & n\alpha + \alpha_3
\end{pmatrix}
\begin{pmatrix}
	[\mu_n] \\[2pt]
	[\lambda_n]
\end{pmatrix}, 
\end{eqnarray*}

and it follows that 

\begin{eqnarray*}
\begin{pmatrix}
	[\mu_n] \\[2pt]
	[\lambda_n]
\end{pmatrix}
=
\begin{pmatrix}
	n\alpha + \alpha_3 & n\beta + \beta_3 \\[2pt]
	-\alpha & -\beta
\end{pmatrix}
\begin{pmatrix}
	[s] \\[2pt]
	[t]
\end{pmatrix}. 
\end{eqnarray*}

Hence, the Seifert invariant of the fiber $c_n$ in $K_n(m_n)$ is $(n\beta+\beta_3) / (n\alpha+\alpha_3)$, 
and $K_n(m_n) = S^2(b, r_1, r_2, (n\beta+\beta_3) / (n\alpha+\alpha_3))$. 

\medskip

\begin{remark}
\label{limit}
Since $(n\beta+\beta_3) / (n\alpha+\alpha_3)$ converges to $\beta/ \alpha$ when $|n|$ tends to $\infty$, 
$M_c(K, m)$ can be regarded as the limit of $K_n(m_n)$ when $|n|$ tends to $\infty$. 
\end{remark}

We divide into three cases: $r_c = -1/0$, $r_c \in \mathbb{Z}$ or $r_c \in \mathbb{Q} \setminus \mathbb{Z}$. 
Except for the last case, 
we do not need the assumption that $M_c(K, m)$ is an L-space. 

\medskip

\noindent
\textbf{Case 1}.\ $r_c = \beta/\alpha = -1/0$. 

Since $\alpha_3 > 0$ and $\alpha \beta_3 - \beta\alpha_3 = 1$, 
we have $\alpha_3 = 1,\ \beta = -1$. 
Hence $K_n(m_n)$ is a Seifert fiber space 
$S^2(b, r_1, r_2, (n\beta+\beta_3) / (n\alpha+\alpha_3)) = S^2(b, r_1, r_2, -n+\beta_3)$, 
which is a lens space for any $n \in \mathbb{Z}$. 
Following Claim~\ref{S2xS1} $S^2(b, r_1, r_2, -n+\beta_3)$ is an L-space 
if $n \ne b+ \beta_3 + r_1 + r_2$. 
Thus $(K_n, m_n)$ is an L-space surgery for all $n \in \mathbb{Z}$ except at most $n = b+ \beta_3 + r_1 + r_2$. 

\bigskip

\noindent
Next suppose that $r_c = \beta/\alpha \ne -1/0$.  
Then the Seifert invariant of $c_n$ is 
$$f(n) = \frac{n\beta+\beta_3}{n\alpha+\alpha_3} 
= \frac{\beta}{\alpha} + \frac{\beta_3 - \frac{\beta}{\alpha}\alpha_3}{n\alpha + \alpha_3} 
= r_c + \frac{\beta_3 - r_c \alpha_3}{n\alpha+\alpha_3}.$$

\noindent
Since $\alpha\beta_3-\beta\alpha_3 = \alpha(\beta_3 - r_c\alpha_3) =1$, 
$\alpha$ and $\beta_3 - r_c\alpha_3$ have the same sign.  

\medskip

\noindent
\textbf{Case 2}.\ $r_c \in \mathbb{Z}$. 
We put $r_c =  p$. 
Then we can write $S^2(b, r_1, r_2, r_c) = S^2(b+ p, r_1, r_2)$.  

\medskip

\noindent
(i) If $b \le -p -3$ or $b \ge -p + 1$, 
then Proposition~\ref{r_c=0}(1) shows that $S^2(b, r_1, r_2, f(n)) = S^2(b+ p, r_1, r_2, f(n)-p)$ is an L-space if  
$-1 \le f(n) - p \le 1$, 
i.e. $p -1 \le f(n) \le p+1$. 
Hence $(K_n, m_n)$ is an L-space for all $n$ but $n \in (x_1, x_2)$, 
where $f(x_1) = p-1$ and $f(x_2) = p+1$; 
see Figure~\ref{seiferter_graph12}(i). 

\medskip

\noindent
(ii) If $b = -p -2$, 
then it follows from Proposition~\ref{r_c=0}(2), 
there is an $\varepsilon > 0$ such that 
$S^2(b, r_1, r_2, f(n)) = S^2(b+ p, r_1, r_2, f(n)-p) = S^2(-2, r_1, r_2, f(n)-p)$ is an L-space 
if  $-1 \le f(n)-p \le \varepsilon$. 
Hence $(K_n, m_n)$ is an L-space except for only finitely many $n \in (x_1, x_2)$, 
where $f(x_1) = p-1$, $f(x_2) = p+\varepsilon$; 
see Figure~\ref{seiferter_graph12}(ii).

\begin{figure}[htbp]
\begin{center}
\includegraphics[width=1.0\linewidth]{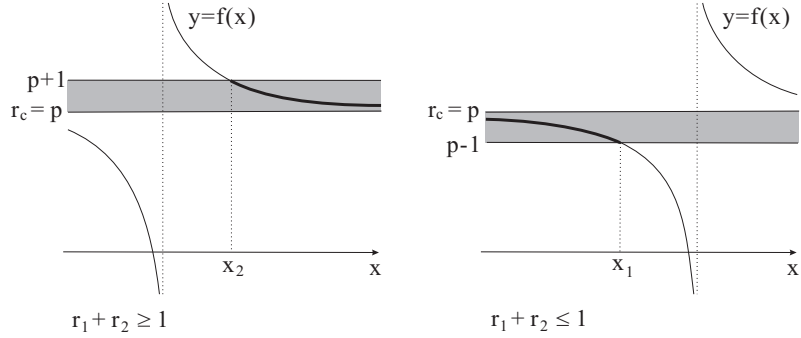}
\caption{$f(x) = \frac{\beta x + \beta_3}{\alpha x + \alpha_3}$}
\label{seiferter_graph12}
\end{center}
\end{figure}

\medskip

 \noindent
(iii) Suppose that $b = -p -1$. 
If $r_1 + r_2 \ge 1$ (resp. $r_1 + r_2 \le 1$), 
then Proposition~\ref{r_c=0}(3) shows that $S^2(b, r_1, r_2, f(n)) = S^2(b+ p, r_1, r_2, f(n)-p) = S^2(-1, r_1, r_2, f(n)-p)$ is an L-space 
if $0 < f(n)-p \le 1$ (resp. $-1 \le f(n)-p < 0$). 
Hence $(K_n, m_n)$ is an L-space for any integer $n \ge x_2$, where $f(x_2) = p+1$ 
(resp. $n \le x_1$, where $f(x_1) = p-1$), 
see Figure~\ref{seiferter_graph3}.

\begin{figure}[htbp]
\begin{center}
\includegraphics[width=0.7\linewidth]{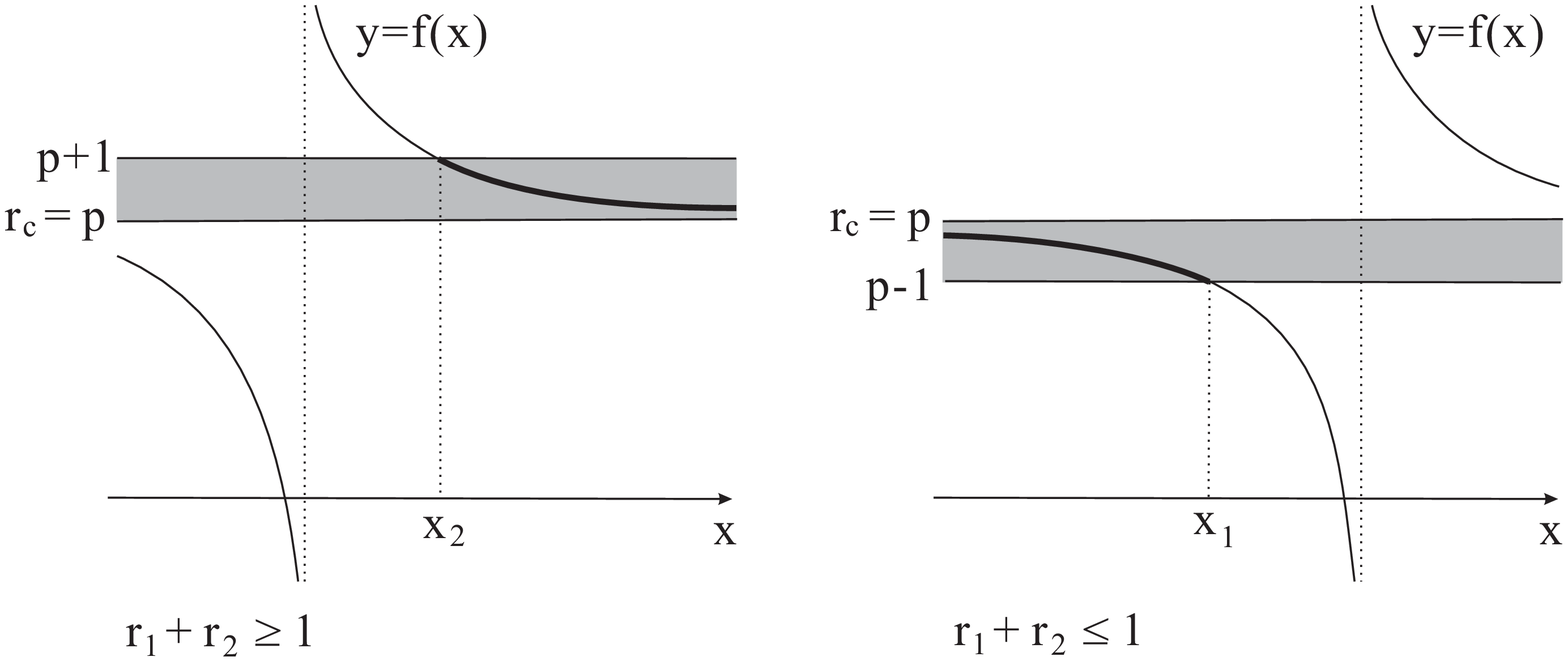}
\caption{$f(x) = \frac{\beta x + \beta_3}{\alpha x + \alpha_3}$}
\label{seiferter_graph3}
\end{center}
\end{figure}

\medskip

\noindent
(iv) If $b = -p$, 
then Proposition~\ref{r_c=0}(4) shows that 
$S^2(b, r_1, r_2, f(n)) = S^2(b+ p, r_1, r_2, f(n)-p) = S^2(r_1, r_2, f(n)-p)$ is an L-space 
if $-\varepsilon \le f(n)-p \le 1$, i.e. 
$p -\varepsilon \le f(n) \le p+1$ for some $\varepsilon > 0$. 
Hence $(K_n, m_n)$ is an L-space for all $n$ but $n \in (x_1, x_2)$, 
where $f(x_1) = p-\varepsilon$ and $f(x_2) = p+1$; 
see Figure~\ref{seiferter_graph12}(iii).

\medskip

\noindent
\textbf{Case 3}.\ $r_c \in \mathbb{Q} \setminus \mathbb{Z}$ and $M_c(K, m) = S^2(b, r_1, r_2, r_c)$ is an L-space. 
We assume $p < r_c  <  p+1$ for some integer $p$.  
Then we have 
$S^2(b, r_1, r_2, r_c) = S^2(b+ p, r_1, r_2, r_c-p)$, 
where $0 < r_c - p < 1$. 

\medskip

\noindent
(i) If $b \le -p-3$ or $b \ge -p+1$, 
then Proposition~\ref{r_c=0}(1) shows that $S^2(b, r_1, r_2, f(n)) = S^2(b+ p, r_1, r_2, f(n)-p)$ is an L-space if  
$-1 \le f(n) - p \le 1$, 
i.e. $p -1 \le f(n) \le p+1$. 
Hence $(K_n, m_n)$ is an L-space for all $n$ but $n \in (x_1, x_2)$, 
where $f(x_1) = p-1$ and $f(x_2) = p+1$; 
see Figure~\ref{seiferter_graph45}(i).

\medskip

\noindent
(ii) Suppose that $b = -p-1$. 
Since $S^2(b, r_1, r_2, r_c) = S^2(b+ p, r_1, r_2, r_c-p) = S^2(-1, r_1, r_2, r_c - p)$ is an L-space, 
by Proposition~\ref{b=-1, -2}(1), 
$S^2(b, r_1, r_2, f(n)) = S^2(-1, r_1, r_2, f(n) - p)$ is an L-space if 
$r_c-p \le f(n) - p \le 1$ (i.e. $r_c \le f(n) \le p+1$). 
Hence $(K_n, m_n)$ is an L-space for any $n \ge x_2$, 
where $f(x_2) = p+1$; 
see Figure~\ref{seiferter_graph45}(ii). 
(Furthermore, 
if $r_1 + r_2 \ge 1$, 
then by Proposition~\ref{r_c=0}(3)(i),  
$S^2(-1, r_1, r_2, f(n) - p)$ is an L-space provided $0 < f(n) - p \le 1$, 
i.e. $p < f(n) \le p+1$. 
Hence $(K_n, m_n)$ is an L-space surgery for any integer $n$ except for $n \in [x_1, x_2)$, 
where $f(x_1) = p$ and $f(x_2) = p+1$. )

\begin{figure}[htbp]
\begin{center}
\includegraphics[width=0.7\linewidth]{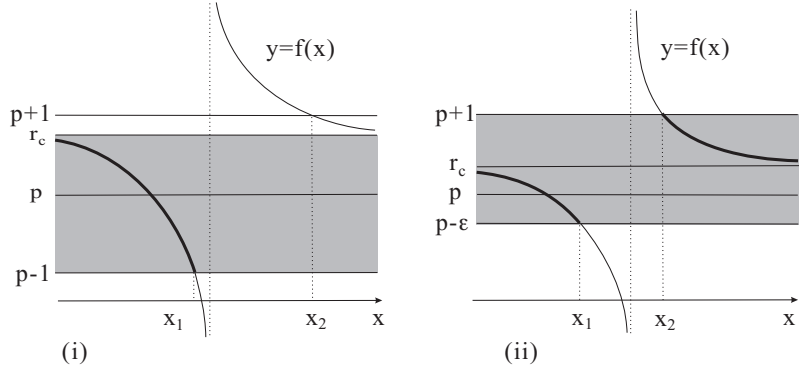}
\caption{$f(x) = \frac{\beta x + \beta_3}{\alpha x + \alpha_3}$}
\label{seiferter_graph45}
\end{center}
\end{figure}

\medskip

\noindent
(iii) Suppose that  $b = -p-2$. 
Since $S^2(b, r_1, r_2, r_c) = S^2(b+ p, r_1, r_2, r_c-p) = S^2(-2, r_1, r_2, r_c - p)$
is an L-space, 
following Proposition~\ref{b=-1, -2}(2), 
$S^2(b, r_1, r_2, f(n)) = S^2(-2, r_1, r_2, f(n) - p)$ is an L-space if 
$-1 \le f(n) - p \le r_c - p$ (i.e. $p-1 \le f(n) \le r_c$). 
Hence $(K_n, m_n)$ is an L-space for any $n \le x_1$, 
where $f(x_1) = p-1$; 
see Figure~\ref{seiferter_graph67}(i).  
(Furthermore, 
if $r_1 + r_2 \le 1$, 
then Proposition~\ref{r_c=0}(2) shows that 
$S^2(-2, r_1, r_2, f(n) - p)$ is an L-space provided $-1 \le f(n) - p < 1$, 
i.e. $p-1 \le f(n) < p+1$. 
Hence $(K_n, m_n)$ is an L-space surgery for any integer $n$ except for $n \in (x_1, x_2]$, 
where $f(x_1) = p-1$ and $f(x_2) = p+1$. )

\medskip

\noindent
(iv) If $b = -p$, 
then Proposition~\ref{r_c=0}(4) shows that 
$S^2(b, r_1, r_2, f(n)) = S^2(b+ p, r_1, r_2, f(n)-p) = S^2(r_1, r_2, f(n)-p)$ is an L-space 
if $-\varepsilon \le f(n)-p \le 1$, i.e. 
$p -\varepsilon \le f(n) \le p+1$ for some $\varepsilon > 0$. 
Hence $(K_n, m_n)$ is an L-space for all $n$ but $n \in (x_1, x_2)$, 
where $f(x_1) = p-\varepsilon$ and $f(x_2) = p+1$; 
see Figure~\ref{seiferter_graph67}(ii). 

\begin{figure}[htbp]
\begin{center}
\includegraphics[width=0.7\linewidth]{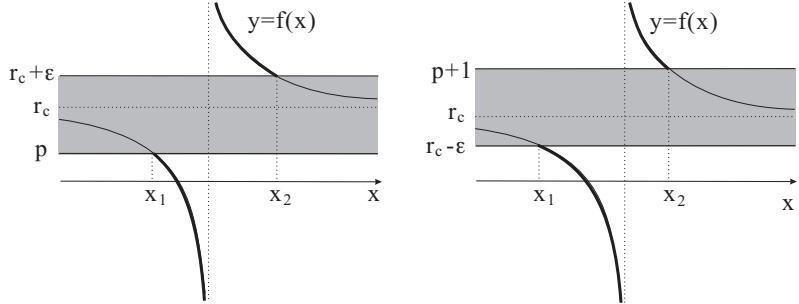}
\caption{$f(x) = \frac{\beta x + \beta_3}{\alpha x + \alpha_3}$}
\label{seiferter_graph67}
\end{center}
\end{figure}

\medskip

Now let us prove the ``only if" part of Theorem~\ref{tree small}. 
We begin by observing: 

\begin{lemma}
\label{McLspace}
$M_c(K, m)$ cannot be $S^2 \times S^1$, 
in particular, 
if $M_c(K, m)$ is a lens space, 
then it is an L-space. 
\end{lemma}

\noindent
\textit{Proof of Lemma~\ref{McLspace}.}
Let $w$ be the linking number between $c$ and $K$.  
Then $H_1(M_c(K, m)) = \langle\ \mu_c,\ \mu_K \ |\ w\mu_c + m\mu_K = 0,\ w\mu_K = 0\ \rangle$, 
where $\mu_c$ is a meridian of $c$ and $\mu_K$ is that of $K$. 
If $M_c(K, m) \cong S^2 \times S^1$, 
then $H_1(M_c(K, m)) \cong \mathbb{Z}$, and we have $w = 0$. 
Let us put $V = S^3 - \mathrm{int}N(c)$, 
which is a solid torus containing $K$ in its interior; $K$ is not contained in any $3$--ball in $V$. 
Since $w = 0$, $K$ is null-homologous in $V$. 
Furthermore, since $c$ is a seiferter for $(K, m)$, 
the result $V(K; m)$ of $V$ after $m$--surgery on $K$ has a (possibly degenerate) Seifert fibration. 
Then \cite[Lemma~3.22]{DMM1} shows that the Seifert fibration of $V(K; m)$ is non-degenerate and 
neither a meridian nor a longitude of $V$ is a fiber in $V(K; m)$, and 
the base surface of $V(K; m)$ is not a M\"obius band.  
Since $K$ is null-homologous in $V$, 
$V(K; m)$ is not a solid torus \cite[Theorem~1.1]{Gabai_solidtorus}, 
and hence, $V(K; m)$ has a Seifert fibration over the disk with at least two exceptional fibers. 
Then $M_c(K, m) = V(K; m) \cup N(c)$ is obtained by attaching $N(c)$ to $V(K; m)$ so that 
the meridian of $N(c)$ is identified with a meridian of $V$. 
Since a regular fiber on $\partial V(K; m)$ intersects a meridian of $V$, 
i.e. a meridian of $N(c)$ more than once, 
$M_c(K, m)$ is a Seifert fiber space over $S^2$ with at least three exceptional fibers. 
Therefore $M_c(K, m)$ cannot be $S^2 \times S^1$. 
This completes a proof. 
\hspace*{\fill} $\square$(Lemma~\ref{McLspace})

\medskip

Suppose first that $K(m)$ is a lens space and $c$ is a core of a genus one Heegaard splitting of $K(m)$. 
Then $V(K; m) = K(m) - \mathrm{int}N(c)$ is a solid torus and 
$M_c(K, m) = V(K; m) \cup N(c)$ is obviously a lens space. 
By Lemma~\ref{McLspace} $M_c(K, m)$ is an L-space. 

In the remaining case, 
as in the proof of the ``if" part of Theorem~\ref{tree small}, 
$M_c(K, m)$ has a form  $S^2(b, r_1, r_2, r_c)$ $(0 < r_1 \le r_2 < 1)$. 

\begin{claim}
\label{rc}
If $r_c = -1/0$ or $r_c \in \mathbb{Z}$, 
then $M_c(K, m)$ is an L-space. 
\end{claim}

\noindent
\textit{Proof of Claim~\ref{rc}.}
If $r_c = -1/0$, 
then $M_c(K, m) = S^2(b, r_1, r_2, -1/0)$ is a connected sum of two lens spaces. 
Since a connected sum of L-spaces is also an L-space \cite[8.1(5)]{Szabo} (\cite{OS2}), 
$M_c(K, m)$ is an L-space. 
If $r_c \in \mathbb{Z}$, 
then $M_c(K, m)$ is a lens space, 
hence it is an L-space by Lemma~\ref{McLspace}. 
\hspace*{\fill} $\square$(Claim~\ref{rc})

\medskip

Now suppose that $M_c(K, m)$ is not an L-space. 
Then by Claim~\ref{rc} $r_c \in \mathbb{Q} \setminus \mathbb{Z}$. 
We write $r_c = r'_c + p$ so that $0 < r'_c < 1$ and $p \in \mathbb{Z}$. 
Then $M_c(K, m) = S^2(b, r_1, r_2, r_c) = S^2(b+p, r_1, r_2, r'_c)$. 
Since $M_c(K, m)$ is not an L-space, 
$b+p =-1$ or $-2$ (Theorem~\ref{SeifertLspace}). 
It follows from Proposition~\ref{nonLspace} that there is an $\varepsilon > 0$ such that 
$K_n(m_n) = S^2(b, r_1, r_2, f(n)) = S^2(b+p, r_1, r_2, f(n)-p) = S^2(-1, r_1, r_2, f(n) -p)$ 
(resp. $S^2(-2, r_1, r_2, f(n) -p)$) 
is not an L-space if $0 < f(n) - p < r'_c + \varepsilon$, i.e. $p < f(n) < r_c + \varepsilon$ 
(resp. $r'_c - \varepsilon < f(n) - p < 1$, i.e. $r_c - \varepsilon < f(n) < p+1$). 
Hence there are at most finitely many integers $n$ such that $K_n(m_n)$ is an L-space, 
i.e. $(K_n, m_n)$ is an L-space surgery. 
See Figure~\ref{nonLspace_graph}. 

\begin{figure}[htbp]
\begin{center}
\includegraphics[width=0.7\linewidth]{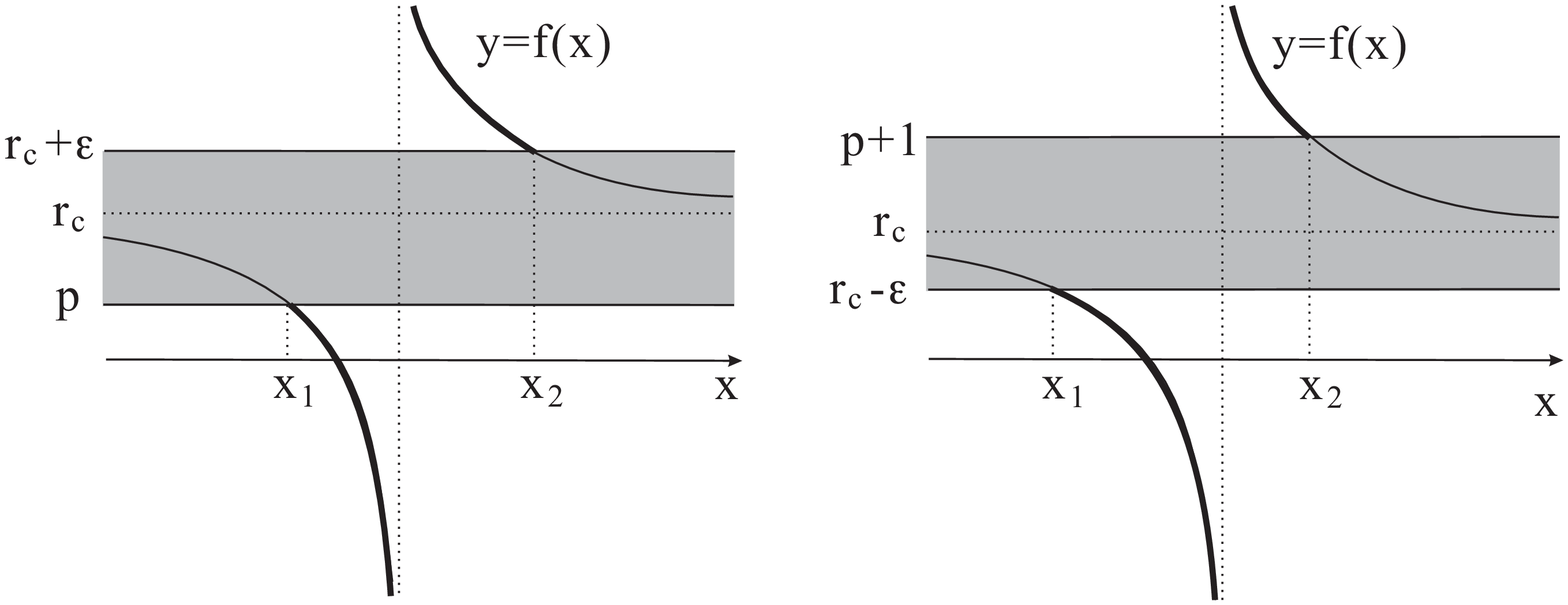}
\caption{$f(x) = \frac{\beta x + \beta_3}{\alpha x + \alpha_3}$}
\label{nonLspace_graph}
\end{center}
\end{figure}

This completes a proof of Theorem~\ref{tree small}. 
\hspace*{\fill} $\square$(Theorem~\ref{tree small})

\bigskip

\noindent
\textit{Proof of Theorem~\ref{tree projective}.}
Note that $K_n(m_n)$ is a Seifert fiber space which admits a Seifert fibration over $\mathbb{R}P^2$, 
or $K_n(m_n)$ has $S^2 \times S^1$ as a connected summand according as $c$ becomes a non-degenerate fiber, 
or a degenerate fiber in $K_n(m_n)$, respectively. 
In the former case, Boyer, Gordon and Watson \cite[Proposition~5]{BGW} prove that $K_n(m_n)$ is an L-space. 
In the latter case, 
$(K_n, m_n)=(O, 0)$ \cite[Theorem~8.1]{GabaiIII}, 
which is not an L-space surgery, 
but there is at most one such integer $n$ \cite[Theorem~5.1]{DMM1}. 
This completes a proof. 
\hspace*{\fill} $\square$(Theorem~\ref{tree projective})

\medskip

\begin{example}
\label{seiferter_projective}
Let us consider a three component link $O \cup c_1 \cup c_2$ depicted in Figure~\ref{projective}. 
It is shown in \cite[Lemma~9.26]{DMM1} that $c_1, c_2$ become fibers in a Seifert fibration of $O(0)$.  
Let $A$ be an annulus in $S^3$ cobounded by $c_1$ and $c_2$. 
Performing $(-l)$--annulus twist along $A$, 
equivalently performing $(1/l+3)$--, $(-1/l+ 3)$--surgeries on $c_1$, $c_2$ respectively, 
we obtain a knot $K_l$ given by Eudave-Mu\~noz \cite{EM2}.  
Then, as shown in \cite{EM2}, 
$(K_l, 12l^2 -4l)$ is a Seifert surgery such that $K_l(12l^2 -4l)$ is a Seifert fiber space over $\mathbb{R}P^2$ 
with at most two exceptional fibers $c_1$, $c_2$ of indices $| l |$, $| -3l+1|$ for $l \ne 0$, 
where we use the same symbol $c_i$ to denote the image of $c_i$ after $(-l)$--annulus twist along $A$. 
Let $c$ be one of $c_1$ or $c_2$.  
Then $c$ is a seiferter for $(K_l, 12l^2 -4l)$.  
Theorem~\ref{tree projective} shows that a knot $K_{l, n}$ obtained from $K_l$ by $n$--twist along $c$ is an L-space knot for all integers $n$. 

\begin{figure}[htbp]
\begin{center}
\includegraphics[width=0.3\linewidth]{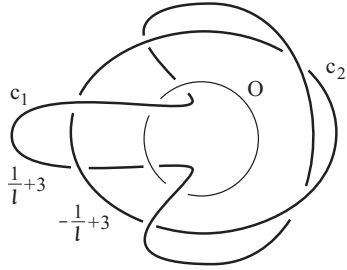}
\caption{$c_1$ and $c_2$ become fibers in a Seifert fibration of $O(0)$.}
\label{projective}
\end{center}
\end{figure}

\end{example}

\bigskip

\section{L-space surgeries and twisting along seiferters II -- degenerate case}
\label{Twisting L-space seiferter II}

In this section we will prove Theorem~\ref{tree degenerate}. 

\medskip

\noindent
\textit{Proof of Theorem~\ref{tree degenerate}.}
Since $K(m)$ has a degenerate Seifert fibration, 
it is a lens space or a connected sum of two lens spaces \cite[Proposition~2.8]{DMM1}. 

\medskip

\noindent
$(1)$\ $K(m)$ is a lens space with degenerate Seifert fibration. 

Then there are at most two degenerate fibers in $K(m)$ \cite[Proposition~2.8]{DMM1}. 
Assume that there are exactly two degenerate fibers. 
Then $(K, m) = (O, 0)$ and the exterior of these two degenerate fibers is $S^1 \times S^1 \times [0, 1]$.  
If $c$ is a non-degenerate fiber, 
then $K_n(m_n)$ has $S^2 \times S^1$ as a connected summand for all integers $n$, 
and thus $(K_n, m_n) = (O, 0)$ for all integers $n$ \cite[Theorem~8.1]{GabaiIII}. 
This contradicts \cite[Theorem~5.1]{DMM1}. 
If $c$ is one of the degenerate fibers, 
then $(K_n, m_n)$ is a lens space, 
which is $S^2 \times S^1$ only when $(K_n, m_n) = (O, 0) = (K_0, m_0)$, 
i.e. $n=0$ \cite[Theorem~5.1]{DMM1}. 
Thus $(K_n, m_n)$ is an L-space surgery except when $n=0$. 

Suppose that $K(m)$ has exactly one degenerate fiber $t_d$. 
There are two cases to consider: 
$K(m) - \mathrm{int}N(t_d)$ is a fibered solid torus or has a non-degenerate Seifert fibration over the M\"obius band with 
no exceptional fiber (\cite[Proposition~2.8]{DMM1}).  
In either case, a meridian of $t_d$ is identified with a regular fiber on $\partial(K(m) - \mathrm{int}N(t_d))$. 

Assume that $K(m) - \mathrm{int}N(t_d)$ is a fibered solid torus. 
Suppose that $c$ is a non-degenerate fiber. 
If $c$ is a core of the solid torus, then $K(m) - \mathrm{int}N(c)$ is a solid torus and $K_n(m_n)$ is a lens space. 
Hence $(K_n, m_n)$ is an L-space surgery except when $K_n(m_n) \cong S^2 \times S^1$, i.e. $(K_n, m_n) = (O, 0)$. 
By \cite[Theorem~5.1]{DMM1} there is at most one such integer $n$. 
If $c$ is not a core in the fibered solid torus  $K(m) - \mathrm{int}N(t_d)$, 
then $K_n(m_n)$ is a lens space ($\not\cong S^2 \times S^1$), a connected sum of two lens spaces, 
or a connected sum of $S^2 \times S^1$ and a lens space ($\not\cong S^3,\ S^2 \times S^1$). 
The last case cannot happen for homological reasons, 
and hence $(K_n, m_n)$ is an L-space surgery. 
If $c$ is the degenerate fiber $t_d$, 
then $K_n(m_n)$ is a lens space, and except for at most integer $n_0$ with $(K_{n_0}, m_{n_0}) = (O, 0)$, 
$(K_n, m_n)$ is an L-space surgery. 

Next consider the case where $K(m) - \mathrm{int}N(t_d)$ has a non-degenerate Seifert fibration over the M\"obius band. 
Then $(K, m) = (O, 0)$; see \cite[Proposition~2.8]{DMM1}. 
If $c$ is a non-degenerate fiber,  
$K_n(m_n)$ has $S^2 \times S^1$ as a connected summand for all integers $n$. 
This implies that $(K_n, m_n) = (O, 0)$ for all $n$ \cite[Theorem~8.1]{GabaiIII}, 
contradicting \cite[Theorem~5.1]{DMM1}.  
Thus $c$ is a degenerate fiber, 
and $K_n(m_n)$ $(n \ne 0)$ is a Seifert fiber space over $\mathbb{R}P^2$ with at most one exceptional fiber, which has finite fundamental group. 
Hence for any non-zero integer $n$, 
$(K_n, m_n)$ is an L-space \cite[Proposition~2.3]{OS3}.  
It follows that if $c$ is a fiber in a degenerate Seifert fibration of a lens space $K(m)$, 
then $(K, m)$ is an L-space surgery except for at most one integer $n$.

\medskip

\noindent
$(2)$\ $K(m)$ is a connected sum of two lens spaces. 

It follows from \cite[Proposition~2.8]{DMM1} that $K(m)$ has exactly one degenerate fiber $t_d$ and $K(m) - \mathrm{int}N(t_d)$ 
is a Seifert fiber space over the disk with two exceptional fibers.  
Note that a meridian of $t_d$ is identified with a regular fiber on $\partial (K(m) - \mathrm{int}N(t_d))$. 
We divide into two cases according as $c$ is a non-degenerate fiber or a degenerate fiber. 

\medskip 
\noindent
(i) $c$ is a non-degenerate fiber.

By \cite[Corollary 3.21(1)]{DMM1} $c$ is not a regular fiber. 
Hence $c$ is an exceptional fiber, 
and $K_n(m_n)$ is a lens space $(\not\cong S^2 \times S^1)$, 
a connected sum of two lens spaces, 
or a connected sum of $S^2 \times S^1$ and a lens space ($\not\cong S^3,\ S^2 \times S^1$). 
The last case cannot happen for homological reasons.  
Hence $(K_n, m_n)$ is an L-space surgery for any integer $n$.  

\medskip

\noindent
(ii) $c$ is a degenerate fiber, i.e. $c = t_d$. 

As in the proof of Theorem~\ref{tree small}, 
let $E$ be $K(m) - \mathrm{int}N(c)$ with a fibered tubular neighborhood of the union of two exceptional fibers $t_1, t_2$
and one regular fiber $t_0$ removed. 
Then $E$ is a product circle bundle over the fourth punctured sphere. 
Take a cross section of $E$ 
such that $K(m)$ has a Seifert invariant 
$S^2(b, r_1, r_2, 1/0)$, 
where the Seifert invariant of $t_0$ is $b \in \mathbb{Z}$, 
that of $t_i$ is $ 0 < r_i < 1$ $(i = 1, 2)$, 
and that of $c$ is $1/0$. 
We may assume that $r_1 \le r_2$. 
Let $s$ be the boundary curve on $\partial N(c)$ of 
the cross section so that $[s]\cdot[t] = 1$ for a regular fiber $t \subset \partial N(c)$. 
Then $[\mu] = [t] \in H_1(\partial N(c))$ and $[\lambda] = -[s] -\beta[t]  \in H_1(\partial N(c))$ 
for some integer $\beta$, i.e. we have: 
\begin{eqnarray*}
\begin{pmatrix}
	[\mu] \\[2pt]
	[\lambda]
\end{pmatrix}
=
\begin{pmatrix}
	0 & 1 \\[2pt]
	-1 & -\beta
\end{pmatrix}
\begin{pmatrix}
	[s] \\[2pt]
	[t]
\end{pmatrix}
\end{eqnarray*} 

Let $c_n$ be the image of $c$ after $n$--twist along $c$. 
Then the argument in the proof of Theorem~\ref{tree small} shows that 
a preferred meridian-longitude pair $(\mu_n, \lambda_n)$ of $\partial N(c_n)$ has the expression: 

\begin{eqnarray*}
\begin{pmatrix}
	[\mu_n] \\[2pt]
	[\lambda_n]
\end{pmatrix}
=
\begin{pmatrix}
	n  & n\beta + 1 \\[2pt]
	-1 & -\beta
\end{pmatrix}
\begin{pmatrix}
	[s] \\[2pt]
	[t]
\end{pmatrix}
\end{eqnarray*}

Thus $K_n(m_n) = S^2(b, r_1, r_2, (n\beta+1)/n)
= S^2(b+\beta, r_1, r_2, (n\beta + 1) /n -\beta) = S^2(b+\beta, r_1, r_2, 1/n)$ for non-zero integer $n$. 

\begin{claim}
\label{K-1K0K1}
$K_n(m_n)$ is an L-space for $n = 0, \pm 1$. 
\end{claim}

\noindent
\textit{Proof of Claim~\ref{K-1K0K1}.}
Recall that $K_0(m_0) = K(m)$ is a connected sum of two lens spaces $L_1$ and $L_2$ 
such that $H_1(L_1) \cong \mathbb{Z}_{\alpha_1}$ and  $H_1(L_2) \cong \mathbb{Z}_{\alpha_2}$, 
where $r_i = \beta_i/\alpha_i$.  
Thus $K_0(m_0)$ is an L-space. 
Since $K_{-1}(m_{-1})$ and $K_1(m_1)$ are lens spaces, 
it remains to show that they are not $S^2 \times S^1$. 
Assume for a contradiction that 
$K_1(m_1)$ or $K_{-1}(m_{-1})$ is $S^2 \times S^1$. 
Then Claim~\ref{S2xS1} shows that 
$r_1 + r_2 = 1$, hence $r_2 = \beta_2/\alpha_2 = (\alpha_1 - \beta_1)/\alpha_1$. 
Thus $\alpha_1 = \alpha_2$, and $H_1(K_0(m_0)) \cong \mathbb{Z}_{\alpha_1} \oplus \mathbb{Z}_{\alpha_2}$ is not cyclic, 
a contradiction. 
Hence neither $K_1(m_1)$ nor $K_{-1}(m_{-1})$ is $S^2 \times S^1$ and 
they are L-spaces. 
\hspace*{\fill} $\square$(Claim~\ref{K-1K0K1})

\medskip

$(1)$ If $b + \beta \le -3$ or $b + \beta \ge 1$, 
Proposition~\ref{r_c=0}(1) shows that 
$K_n(m_n) = S^2(b+\beta, r_1, r_2, 1/n)$ is an L-space 
if $-1 \le 1/n \le 1$, 
i.e. $n \le -1$ or $n \ge 1$. 
See Figure~\ref{lenslens_graph2}(i). 
Since $K_0(m_0)$ is also an L-space (Claim~\ref{K-1K0K1}), 
$K_n(m_n)$ is an L-space for any integer $n$. 

\medskip

$(2)$ If $b + \beta = -2$, 
Proposition~\ref{r_c=0}(2) shows that there is an $\varepsilon > 0$ such that 
$K_n(m_n) = S^2(b+\beta, r_1, r_2, 1/n)$ is an L-space 
if $-1 \le 1/n \le \varepsilon$. 
Hence $K_n(m_n)$ is an L-space 
if $n \le -1$ or $n \ge 1/\varepsilon$.  
See Figure~\ref{lenslens_graph2}(ii). 
This, together with Claim~\ref{K-1K0K1}, shows that 
$K_n(m_n)$ is an L-space 
if $n \le 1$ or $n \ge 1/\varepsilon$. 
\medskip 

$(3)$ Suppose that $b + \beta = -1$. 
Then Proposition~\ref{r_c=0}(3) shows that if $r_1 + r_2 \ge 1$ (resp. $r_1 + r_2 \le 1$), 
$K_n(m_n) = S^2(b+\beta, r_1, r_2, 1/n)$ is an L-space 
for any integer $n$ satisfying $0 < 1/n \le 1$ (resp. $-1 \le 1/n < 0$), 
i.e. $n \ge 1$ (resp. $n \le -1$). 
See Figure~\ref{lenslens_graph2}(i). 
Combining Claim~\ref{K-1K0K1}, 
we see that $K_n(m_n)$ is an L-space for any $n \ge -1$ (resp. $n \le 1$).

\medskip

$(4)$ If $b + \beta = 0$,  
then Proposition~\ref{r_c=0}(4) shows that there is an $\varepsilon > 0$ such that 
$K_n(m_n) = S^2(b+\beta, r_1, r_2, 1/n)$ is an L-space 
if $-\varepsilon \le 1/n \le 1$. 
Hence $K_n(m_n)$ is an L-space 
if $n \ge 1$ or $n \le -1/\varepsilon$. 
See Figure~\ref{lenslens_graph2}(iii). 
This, together with Claim~\ref{K-1K0K1}, shows that 
$K_n(m_n)$ is an L-space 
if $n \ge -1$ or $n \le -1/\varepsilon$. 

\medskip

This completes a proof of Theorem~\ref{tree degenerate}. 
\hspace*{\fill} $\square$(Theorem~\ref{tree degenerate})

\medskip

\begin{figure}[htbp]
\begin{center}
\includegraphics[width=1.0\linewidth]{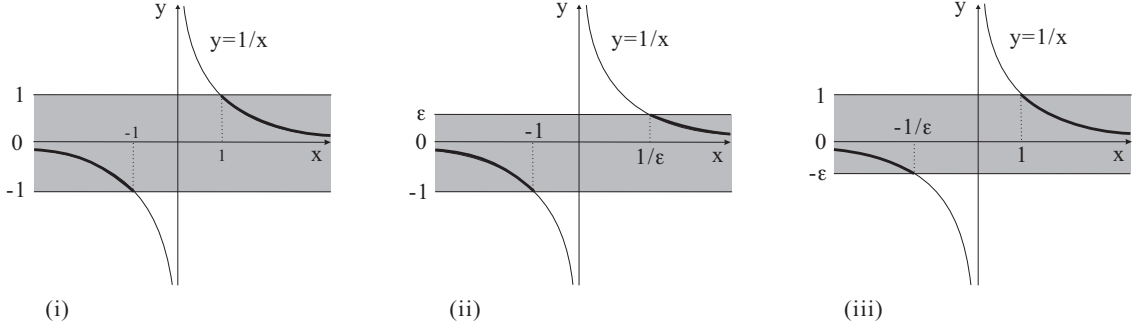}
\caption{(i) $-1 \le 1/n \le 1$ if $n \le -1$ or $n \ge 1$, 
(ii) $-1 \le 1/n \le \varepsilon$ if $n \le -1$ or $n \ge 1/\varepsilon$, 
(iii) $-\varepsilon \le 1/n \le 1$ if $n \le -1/\varepsilon$ or $n \ge 1$.}
\label{lenslens_graph2}
\end{center}
\end{figure}

As shown by Greene \cite[Theorem~1.5]{Greene}, 
if $K(m)$ is a connected sum of lens spaces, 
then $K$ is a torus knot or a cable of a torus knot. 
More precisely, 
$(K, m) = (T_{p, q}, pq)$ or $(C_{p, q}(T_{r, s}), pq)$, where $p = qrs \pm 1$. 
Note that $T_{p, q}(pq) = L(p, q) \sharp L(q, p)$ and $C_{p, q}(T_{r, s})(pq) = L(p, qs^2) \sharp L(q, \pm 1)$. 

Let us continue to prove Theorem~\ref{Tpq} which is a refinement of Theorem~\ref{tree degenerate}(2). 

\medskip

\noindent
\textit{Proof of Theorem~\ref{Tpq}.}
In the following $(K, m)$ is either $(T_{p, q}, pq)$ or $(C_{p, q}(T_{r, s}), pq)$, 
where $p, q \ge 2$ and $p = qrs \pm 1$. 
If $c$ becomes a non-degenerate fiber in $K(m)$, 
then as shown in the proof of Theorem~\ref{tree degenerate}, 
$K_n$ is an L-space knot for any integer $n$. 
So we assume that $c$ becomes a degenerate fiber in $K(m)$. 
Recall from Theorem~3.19(3) in \cite{DMM1} that the linking number $l$ between $c$ and $K$ is not zero.  
Recall also that $K_n(m_n)$ is expressed as 
$S^2(b+\beta, r_1, r_2, 1/n) = S^2(b+\beta, \beta_1/\alpha_1, \beta_2/\alpha_2, 1/n)$, 
where $0 < r_i = \beta_i / \alpha_i < 1$ and $\alpha_i \ge 2$. 
See the proof of Theorem~\ref{tree degenerate}. 
Note that $\{\alpha_1, \alpha_2 \} = \{ p, q \}$, 
and $\alpha_1\alpha_2 = pq \ge 6$.

\begin{claim}
\label{b+beta not -2}
$b + \beta \ne -2$. 
\end{claim}

\noindent
\textit{Proof of Claim~\ref{b+beta not -2}.}
Assume for a contradiction that $b + \beta = -2$. 
Then $K_1(m_1) = S^2(-2,\ \beta_1/\alpha_1,\ \beta_2/\alpha_2,\ 1) = S^2(-1,\ \beta_1/\alpha_1,\ \beta_2/\alpha_2)$. 
Hence $| H_1(K_1(m_1))| = |-\alpha_1\alpha_2 + \alpha_1\beta_2 + \alpha_2\beta_1|$, 
which coincides with $pq + l^2 = \alpha_1\alpha_2 + l^2$. 
Since $\alpha_1\alpha_2 + l^2 > \alpha_1\alpha_2$, 
we have  $|-\alpha_1\alpha_2 + \alpha_1\beta_2 + \alpha_2\beta_1| > \alpha_1\alpha_2$. 
This then implies 
$\beta_1 / \alpha_1 + \beta_2 / \alpha_2 > 2$ or 
$\beta_1 / \alpha_1 + \beta_2 / \alpha_2 < 0$. 
Either case cannot happen, because $0 < \beta_i / \alpha_i < 1$. 
Thus $b + \beta \ne -2$.
\hspace*{\fill} $\square$(Claim~\ref{b+beta not -2})

\medskip

\begin{claim}
\label{b+beta-1}
If $b + \beta = -1$, $\beta_1/\alpha_1 + \beta_2 / \alpha_2 > 1$. 
\end{claim}

\noindent
\textit{Proof of Claim~\ref{b+beta-1}.}
If $b + \beta = -1$, 
then $K_1(m_1) = S^2(-1, \beta_1/\alpha_1, \beta_2/\alpha_2, 1) = S^2(\beta_1/\alpha_1, \beta_2/\alpha_2)$.  
Thus $|H_1(K_1(m_1))| = \alpha_1\beta_2 + \alpha_2\beta_1$, 
which coincides with $pq + l^2 = \alpha_1\alpha_2 + l^2$.  
Since $\alpha_1\alpha_2 + l^2 > \alpha_1\alpha_2$, 
we have $\alpha_1\beta_2 + \alpha_2\beta_1 > \alpha_1\alpha_2$. 
This shows $\beta_1 / \alpha_1 + \beta_2 / \alpha_2 > 1$. 
\hspace*{\fill} $\square$(Claim~\ref{b+beta-1})

\medskip

Claims~\ref{b+beta not -2} and \ref{b+beta-1}, 
together with the argument in the proof of Theorem~\ref{tree degenerate} prove 
that $K_n$ is an L-space knot for any $n \ge -1$. 

\medskip

Now let us prove that $K_n$ is an L-space knot for all integers $n$ under the assumption $l^2 \ge  2pq$. 

\begin{claim}
\label{b+beta not -1}
If $l^2 \ge 2pq$, 
then $b + \beta \ne -1$. 
\end{claim}

\noindent
\textit{Proof of Claim~\ref{b+beta not -1}.}
Assume that $l^2 \ge 2pq$, 
and suppose for a contradiction $b + \beta = -1$. 
Then $K_{-1}(m_{-1}) = S^2(-1, \beta_1/\alpha_1, \beta_2/\alpha_2, -1) = S^2(-2,  \beta_1/\alpha_1, \beta_2/\alpha_2)$, 
and $|H_1(K_{-1}(m_{-1}) )| = |-2\alpha_1\alpha_2 + \alpha_1\beta_2 + \alpha_2\beta_1|$, 
which coincides with $|pq - l^2|$.  
The assumption $l^2 \ge 2pq = 2\alpha_1\alpha_2$ implies that  
$|pq-l^2| = l^2 - pq = l^2 - \alpha_1\alpha_2 \ge \alpha_1\alpha_2$.  
Hence $|-2\alpha_1\alpha_2 + \alpha_1\beta_2 + \alpha_2\beta_1|  = |pq - l^2| \ge \alpha_1\alpha_2$. 
Thus we have $\beta_1/\alpha_1 + \beta_2 / \alpha_2 \ge 3$ or $\beta_1/\alpha_1 + \beta_2 / \alpha_2 \le 1$. 
The former case cannot happen because $0 < \beta_i / \alpha_i < 1$, 
and the latter case contradicts Claim~\ref{b+beta-1} which asserts $\beta_1/\alpha_1 + \beta_2 / \alpha_2 > 1$.  
Hence $b + \beta \ne -1$. 
\hspace*{\fill} $\square$(Claim~\ref{b+beta not -1})

\medskip

\begin{claim}
\label{b+beta not 0}
If $l^2 \ge 2pq$, 
then $b + \beta \ne 0$. 
\end{claim}

\noindent
\textit{Proof of Claim~\ref{b+beta not 0}.}
Suppose for a contradiction that 
$b + \beta = 0$. 
Then $K_{-1}(m_{-1}) = S^2(0, \beta_1/\alpha_1, \beta_2/\alpha_2, -1) = S^2(-1, \beta_1/\alpha_1, \beta_2/\alpha_2)$, 
and $|H_1(K_{-1}(m_{-1}) )| = |-\alpha_1\alpha_2 + \alpha_1\beta_2 + \alpha_2\beta_1|$, 
which coincides with $|pq - l^2|$. 
Since $l^2 \ge 2pq = 2\alpha_1\alpha_2$, 
$|pq - l^2| = l^2 - pq = l^2 - \alpha_1\alpha_2 \ge \alpha_1\alpha_2$. 
Thus we have 
$|-\alpha_1\alpha_2 + \alpha_1\beta_2 + \alpha_2\beta_1| = |pq-l^2| \ge \alpha_1\alpha_2$.  
This then implies 
$\beta_1 / \alpha_1 + \beta_2 / \alpha_2 \ge 2$ or 
$\beta_1 / \alpha_1 + \beta_2 / \alpha_2 \le 0$. 
Either case cannot happen, because $0 < \beta_i / \alpha_i < 1$. 
Thus $b + \beta \ne 0$. 
\hspace*{\fill} $\square$(Claim~\ref{b+beta not 0})

\medskip

Under the assumption $l^2 \ge 2pq$, 
Claims~\ref{b+beta not -2}, \ref{b+beta not -1} and \ref{b+beta not 0} imply that 
$b+\beta \le -3$ or $b + \beta \ge 1$. 
Then the proof of Theorem~\ref{tree degenerate} 
enables us to conclude that $K_n$ is an L-space knot for all integers $n$. 
\hspace*{\fill} $\square$(Theorem~\ref{Tpq})

\medskip

\begin{example}
\label{trefoil_c6}
Let $K$ be a torus knot $T_{3, 2}$ and $c$ an unknotted circle depicted in Figure~\ref{T32_5}; 
the linking number between $c$ and $T_{3, 2}$ is $5$. 
Then $c$ coincides with $c^+_{3, 2}$ in Section~\ref{twisted torus knots}, 
and it is a seiferter for $(T_{3, 2}, 6)$. 
Let $K_n$ be a knot obtained from $T_{3, 2}$ by $n$--twist along $c$. 
Since $5^2 > 2\cdot 3 \cdot 2 = 12$, 
following Theorem~\ref{Tpq} $K_n$ is an  L-space knot for all integers $n$. 
\end{example}

\begin{figure}[htbp]
\begin{center}
\includegraphics[width=0.25\linewidth]{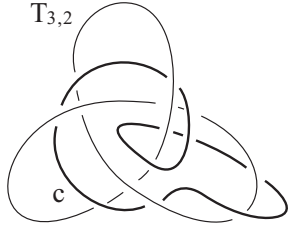}
\caption{$c$ is a seiferter for $(T_{3, 2}, 6)$.}
\label{T32_5}
\end{center}
\end{figure}

\medskip

Example~\ref{cable} below gives an example of a seiferter for $(K, m)$, 
where $K$ is a cable of a torus knot and $K(m)$ is a connected sum of two lens spaces. 

\begin{example}
\label{cable}
Let $k$ be a Berge knot $Spor{\bf a}[p]$ $(p > 1)$. 
Then $k(22p^2+9p+1)$ is a lens space and \cite[Proposition~8.1 and Table~9]{DMM2} shows that 
$(k, 22p^2+9p+1)$ has a seiferter $c$ such that the linking number between $c$ and $k$ is $4p+1$ 
and $(-1)$--twist along $c$ converts $(k, 22p^2+9p+1)$ into 
$(C_{6p+1, p}(T_{3, 2}), p(6p+1))$. 
Since $p > 1$, 
$C_{6p+1, p}(T_{3, 2})$ is a nontrivial cable of $T_{3, 2}$. 
Thus $c$ is a seiferter for $(C_{6p+1, p}(T_{3, 2}), p(6p+1))$. 
Let $K_n$ be a knot obtained from $C_{6p+1, p}(T_{3, 2})$ by $n$--twist along $c$ so that $K_1 = k$. 
Since $(4p+1)^2 \ge 2(6p+1)p$,  
Theorem~\ref{Tpq} shows that $K_n$ is an L-space knot for all integers $n$. 

Finally we show that  $K_n$ is hyperbolic if $|n| > 3$. 
As shown in \cite[Figure~41]{DMM2}, $K_n$ admits a Seifert surgery yielding a small Seifert space which is not a lens space, 
so we see that $c$ becomes a degenerate fiber in $C_{6p+1, p}(T_{3, 2})(p(6p+1))$ \cite[Lemma~5.6(1)]{DMM1}. 
Hence Corollary~3.21(3) in \cite{DMM1} shows that the link $C_{6p+1, p}(T_{3, 2}) \cup c$ is hyperbolic. 
Now the result follows from \cite[Proposition~5.11(3)]{DMM1}. 
\end{example}

\medskip

We close this section with the following observation, 
which shows the non-uniqueness of degenerate Seifert fibration of 
a connected sum of two lens spaces. 

Let $c$ be a seiferter for $(T_{p, q}, pq)$ which becomes a degenerate fiber in $T_{p, q}(pq)$. 
As the simplest example of such a seiferter $c$, 
take a meridian $c_{\mu}$ of $T_{p, q}$. 
Then $c_{\mu}$ is isotopic to the core of the filled solid torus $($i.e. the dual knot of $T_{p, q}$$)$ in $T_{p, q}(pq)$, 
which is a degenerate fiber. 
Hence $c_{\mu}$ is a seiferter for $(T_{p, q}, pq)$ which becomes a degenerate fiber in $T_{p, q}(pq)$, 
and $T_{p, q} - \mathrm{int}N(c_{\mu})$ is homeomorphic to $S^3 - \mathrm{int}N(T_{p, q})$. 
However, in general, 
$T_{p, q} (pq)- \mathrm{int}N(c)$ is not necessarily homeomorphic to $S^3 - \mathrm{int}N(T_{p, q})$.

\medskip
\begin{example}
\label{not torus knot exterior}
Let us take an unknotted circle $c$ as in Figure~\ref{T53_2_nontorus}.  
Then $c$ is a seiferter for $(T_{5, 3}, 15)$ which becomes a degenerate fiber in $T_{5, 3}(15)$, 
but $T_{5, 3}(15) - \mathrm{int}N(c)$ is not homeomorphic to $S^3 - \mathrm{int}N(T_{5, 3})$. 
\end{example}

\noindent
\textit{Proof of Example~\ref{not torus knot exterior}.}
As shown in Figure~\ref{T53_2_nontorus}, 
$T_{5, 3}(15)$ is the the two-fold branched cover of 
$S^3$ branched along $L'$ and $c$ is the preimage of an arc $\tau$. 
Hence $T_{5, 3}(15) - \mathrm{int}N(c)$ is a Seifert fiber space $D^2(2/3, -2/5)$. 
Since $|H_1(S^2(2/3, -2/5, x))| = |4+15x|$ cannot be $1$ for any integer $x$, 
the Seifert fiber space $T_{5, 3}(15) - \mathrm{int}N(c)$ cannot be embedded in $S^3$, 
and hence it is not homeomorphic to $S^3 - \mathrm{int}N(T_{5, 3})$. 
(Note that $c$ coincides with $c^-_{5, 3}$ in Section~\ref{twisted torus knots}.) 
\hspace*{\fill} $\square$(Example~\ref{not torus knot exterior})

\begin{figure}[htbp]
\begin{center}
\includegraphics[width=0.85\linewidth]{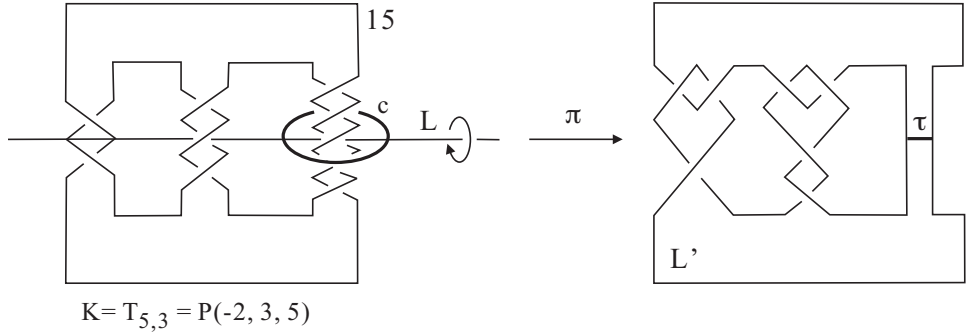}
\caption{$T_{5, 3}(15)$ is the the two-fold branched cover of $S^3$ branched along $L'$.}
\label{T53_2_nontorus}
\end{center}
\end{figure}

\bigskip

\section{L-space twisted torus knots}
\label{twisted torus knots}

Each torus knot has obviously an unknotted circle $c$ which satisfies the 
desired property in Question~\ref{motivation}. 

\begin{example}
\label{obvious example}
Embed a torus knot $T_{p, q}$ into a genus one Heegaard surface of $S^3$. 
Then cores of the Heegaard splitting $s_p$ and $s_q$ are seiferters for $(T_{p, q}, m)$ for all integers $m$. 
We call them \textit{basic seiferters} for $T_{p, q}$; see Figure~\ref{basicseifertersT32}. 
An $n$--twist along $s_p$ (resp. $s_q$) converts $T_{p, q}$ into a torus knot $T_{p+nq, q}$ (resp. $T_{p, q + np}$), 
and hence $n$--twist along a basic seiferter yields an L-space knot for all $n$. 
\end{example}

\begin{figure}[htbp]
\begin{center}
\includegraphics[width=0.23\linewidth]{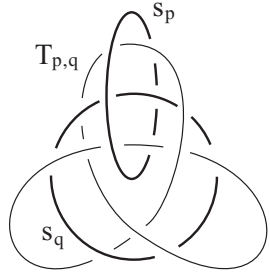}
\caption{$s_p$ and $s_q$ are basic seiferters for $(T_{p, q}, m)$.}
\label{basicseifertersT32}
\end{center}
\end{figure}

Twistings along a basic seiferter keep the property of being L-space knots, 
but produce only torus knots. 
In the following, 
we will give another circle $c$ such that 
twistings $T_{p, q}$ along $c$ produce an infinite family of hyperbolic L-space knots. 

\medskip

\begin{definition}[\textbf{twisted torus knot \cite{D}}]
\label{Dean}
Let $\Sigma$ be a genus one Heegaard surface of $S^3$. 
Let $T_{p, q}$ $(p > q \ge 2)$ be a $(p, q)$--torus knot which lies on $\Sigma$. 
Choose an unknotted circle $c \subset S^3 - T_{p, q}$ so that it bounds a disk $D$ such that 
$D \cap \Sigma$ is a single arc intersecting $T_{p, q}$ in $r$ $(2 \le r \le p+q)$ points in the same direction. 
A \textit{twisted torus knot} $K(p, q; r, n)$ is a knot obtained from $T_{p, q}$ by adding $n$ full twists along $c$. 
\end{definition}

\begin{remark}
\label{p or q}
Twisting $T_{p, q}$ along the basic seiferter $s_p$ $($resp. $s_q$$)$ $n$--times,  
we obtain the twisted torus knot $K(p, q; q, n)$ $($resp. $K(p, q; p, n)$$)$, 
which is a torus knot $T_{p+nq, q}$ $($resp. $T_{p, q+np}$$)$, 
and hence an L-space knot. 
\end{remark}

In \cite{V} Vafaee studies twisted torus knots from a viewpoint of knot Floer homology and 
showed that twisted torus knots $K(p, kp \pm 1; r, n)$, 
where $p \ge 2, k \ge 1, n > 0$ and $0 < r < p$ is an L-space knot if and only if 
either $r = p-1$ or $r \in \{ 2, p-2 \}$ and $n = 1$. 
We will give yet more twisted torus knots which are L-space knots by 
combining seiferter technology and Theorem~\ref{Tpq}. 

\medskip

\noindent
\textit{Proof of Theorem~\ref{twisted torus knot}.} 
In the following, 
let $\Sigma$ be a genus one Heegaard surface of $S^3$, 
which bounds solid tori $V_1$ and $V_2$. 

\medskip

$\bullet$ $K(p, q; p + q, n)$ $(p > q \ge 2)$. \ 
Given any torus knot $T_{p, q}$ $(p > q \ge 2)$ on $\Sigma$, 
let us take an unknotted circle $c^+_{p. q}$ in $S^3 - T_{p, q}$ as depicted in Figure~\ref{cpq_plus}(i); 
the linking number between $c^+_{p, q}$ and $T_{p, q}$ is $p+q$. 

\begin{figure}[h]
\begin{center}
\includegraphics[width=0.75\linewidth]{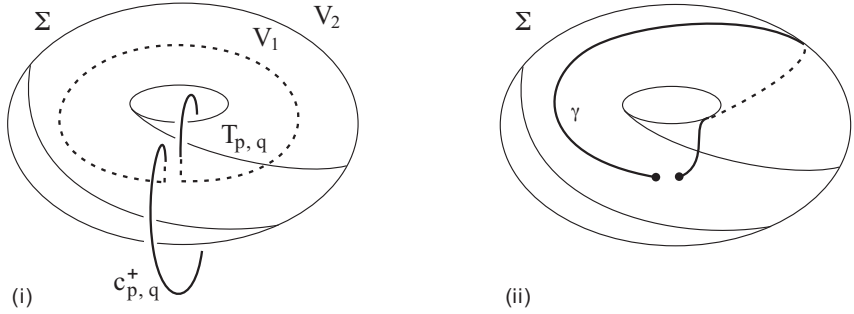}
\caption{$c^+_{p, q}$ is a seiferter for $(T_{p, q}, pq)$.}
\label{cpq_plus}
\end{center}
\end{figure}

Let $V$ be a solid torus $S^3 - \mathrm{int}N(c^+_{p, q})$, 
which contains $T_{p, q}$ in its interior. 
Lemma~9.1 in \cite{MM1} shows that 
$V(K; pq) = T_{p, q}(pq) - \mathrm{int}N(c^+_{p, q})$ is a Seifert fiber space over the disk with two exceptional fibers of indices $p, q$, 
and a meridian of $N(c^+_{p, q})$ coincides with a regular fiber on $\partial V(K; pq)$. 
Hence $c^+_{p, q}$ is a degenerate fiber in $T_{p, q}(pq)$, 
and thus it  is a seiferter for $(T_{p, q}, pq)$. 
Let $D$ be a disk bounded by $c^+_{p, q}$. 
Since the arc $c^+_{p, q} \cap V_i$ is isotoped in $V_i$ to an arc $\gamma \subset \Sigma$ depicted in  
Figure~\ref{cpq_plus}(ii) leaving its endpoints fixed,  
the disk $D$ can be isotoped so that $D \cap \Sigma = \gamma$, 
which intersects $T_{p, q}$ in $p+q$ points in the same direction. 
Thus $n$--twist along $c^+_{p, q}$ converts $T_{p, q}$ into the twisted torus knot $K(p, q; p + q, n)$. 
Since $c^+_{p, q}$ is a seiferter for $(T_{p, q}, pq)$ and $(p+q)^2 = p^2+q^2 + 2pq > 2pq$, 
we can apply Theorem~\ref{Tpq} to conclude that 
$T(p, q, p + q, n)$ is an L-space knots for all integers $n$.  

We show that $T(p, q, p + q, n)$ is hyperbolic if $|n| > 3$. 
By linking number consideration, we see that $c^+_{p, q}$ is not a basic seiferter. 
Then Corollary~3.21(3) in \cite{DMM1} (\cite[Claim~9.2] {MM1}) shows that $T_{p, q} \cup c^+_{p. q}$ is a hyperbolic link. 
Thus \cite[Proposition~5.11(2)]{DMM1} that $K(p, q; p + q, n)$ is a hyperbolic knot if $|n| > 3$. 

\medskip

$\bullet$ $K(p, q; p - q, n)$ $(p > q \ge 2)$. \ 
Suppose that $p-q \ne 1$. 
Then let us take $c^-_{p, q}$ as in Figure~\ref{cpq_minus}(i) instead of $c^+_{p, q}$; 
the linking number between $c^-_{p, q}$ and $T_{p, q}$ is $p-q$. 
It follows from \cite[Remark 4.7]{DMM1} that $c^-_{p, q}$ is also a seiferter for $(T_{p, q}, pq)$ and the link $T_{p, q} \cup c^-_{p, q}$ is hyperbolic. 
Note that if $p-q = 1$, 
then $c^-_{p, q}$ is a meridian of $T_{p, q}$.  
As above we see that each arc $c^-_{p, q} \cap V_i$ is isotoped in $V_i$ to an arc $\gamma \subset \Sigma$ 
depicted in Figure~\ref{cpq_minus}(ii) leaving its endpoints fixed. 
So a disk $D$ bounded by $c^-_{p, q}$ can be isotoped so that $D \cap \Sigma = \gamma$, 
which intersects $T_{p, q}$ in $p-q$ points in the same direction. 
Thus $n$--twist along $c^-_{p, q}$ converts $T_{p, q}$ into the twisted torus knot $K(p, q; p - q, n)$. 
Since $c^-_{p, q}$ is a seiferter for $(T_{p, q}, pq)$,  
Theorem~\ref{Tpq} shows that $T(p, q, p - q, n)$ is an L-space knot for any $n \ge -1$.  
Following \cite[Proposition~5.11(2)]{DMM1} $T(p, q, p - q, n)$ is a hyperbolic knot if $|n| > 3$. 

\begin{figure}[h]
\begin{center}
\includegraphics[width=0.75\linewidth]{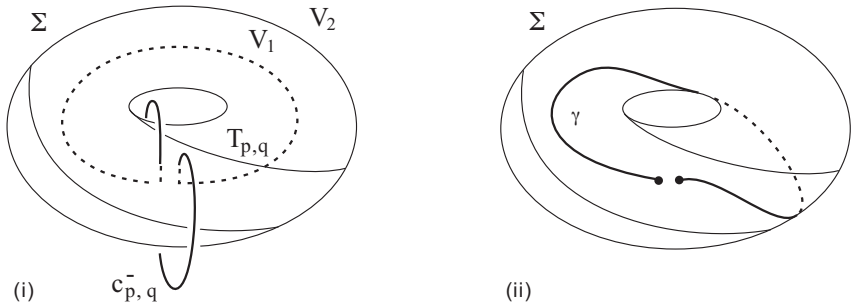}
\caption{$c^-_{p, q}$ is a seiferter for $(T_{p, q}, pq)$.}
\label{cpq_minus}
\end{center}
\end{figure}

\medskip

$\bullet$ $K(3p+1, 2p+1; 4p+1, n)$ $(p > 0)$. \ 
Let $k$ be a torus knot $T_{p, 2p+1}$ on a genus two Heegaard surface, 
with unknotted circles $\alpha$ and $c$ as shown in Figure~\ref{sporbpi1A}. 
Applying $1$--twist along $\alpha$, 
we obtain a torus knot $T_{3p+1, 2p+1}$. 
We continue to use the same symbol $c$ to denote the image of $c$ after $1$--twist along $\alpha$; 
the linking number between $c$ and $T_{3p+1, 2p+1}$ is $4p+1$. 
Note that $1$--twist along $c$ converts $T_{3p+1, 2p+1}$ into a Berge knot $Spor{\bf b}[p]$ as shown in \cite[Subsection 8.2]{DMM2}. 
Following \cite[Lemma~8.4]{DMM2} $c$ is a seiferter for 
a lens space surgery $(Spor{\bf b}[p], 22p^2 + 13p+2) = (Spor{\bf b}[p], (3p+1)( 2p+1) + (4p+1)^2)$. 
Thus $c$ is also a seiferter for $(T_{3p+1, 2p+1}, (3p+1)( 2p+1))$.  
Let $D$ be a disk bounded by $c$. 
Then $T_{3p+1, 2p+1} \cup D$ can be isotoped so that $T_{3p+1, 2p+1}$ lies on $\Sigma$,  
$D \cap \Sigma$ consists of a single arc, 
which intersects $T_{3p+1, 2p+1}$ in $4p+1$ points in the same direction.  
Thus $n$--twist along $c$ converts $T_{3p+1, 2p+1}$ into a twisted torus knot $K(3p+1, 2p+1; 4p+1, n)$. 
Since $c$ is a seiferter for $(T_{3p+1, 2p+1}, (3p+1)(2p+1))$ and $(4p+1)^2 > 2(3p+1)(2p+1)$, 
Theorem~\ref{Tpq} shows that $K(3p+1, 2p+1; 4p+1, n)$ is an L-space knot for all integers $n$. 
Let us observe that $K(3p+1, 2p+1; 4p+1, n)$ is a hyperbolic knot if $|n| > 3$. 
Figure~44 in \cite{DMM2} shows that $n$--twist converts $(T_{3p+1, 2p+1}, (3p+1)(2p+1))$ into a Seifert surgery which is not a lens space surgery if $|n| \ge 2$. 
Hence $c$ becomes a degenerate fiber in $T_{3p+1, 2p+1}((3p+1)(2p+1))$ \cite[Lemma~5.6(1)]{DMM1}, 
and Corollary~3.21(3) in \cite{DMM1} shows that the link $T_{3p+1, 2p+1} \cup c$ is hyperbolic. 
Now the result follows from \cite[Proposition~5.11(2)]{DMM1}. 

\begin{figure}[h]
\begin{center}
\includegraphics[width=0.6\linewidth]{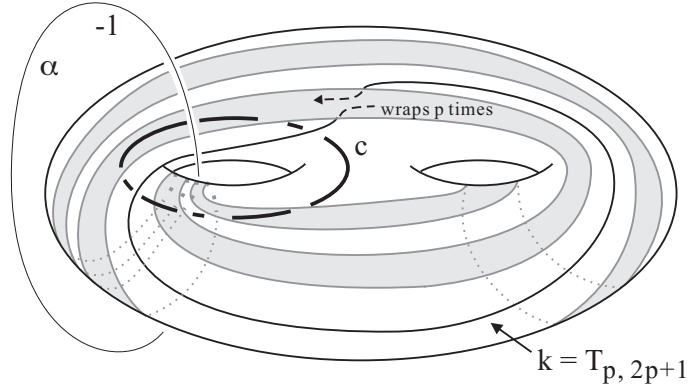}
\caption{A surgery description of $T_{3p+1, 2p+1}$ and a seiferter $c$}
\label{sporbpi1A}
\end{center}
\end{figure}

\medskip

$\bullet$ $K(3p+2, 2p+1; 4p+3, n)$ $(p > 0)$. \ 
As above, we follow the argument in \cite[Subsection~8.3]{DMM2}, but we need to take the mirror image at the end. 
Take a torus knot $k = T_{-p-1, 2p+1}$ on a genus two Heegaard surface of $S^3$,  
unknotted circles $\alpha'$ and $c'$ as shown in Figure~\ref{sporcfiberA}. 
Then $(-1)$--twist along $\alpha'$ converts $T_{-p-1, 2p+1}$ into $T_{-3p-2, 2p+1}$. 
As above we denote the image of $c'$ after $(-1)$--twist along $\alpha'$ by the same symbol $c'$; 
the linking number between $c'$ and $T_{-3p-2, 2p+1}$ is $4p+3$. 
Note that $(-1)$--twist along $c'$ converts $T_{-3p-2, 2p+1}$ into a Berge knot $Spor{\bf c}[p]$ as shown in \cite[Subsection 8.3]{DMM2}. 
Then Lemma~8.6 in \cite{DMM2} shows that $c'$ is a seiferter for 
a lens space surgery $(Spor{\bf c}[p], -22p^2-31p-11) = (Spor{\bf c}[p], (-3p-2)( 2p+1) - (4p+3)^2)$. 
Thus $c'$ is also a seiferter for $(T_{-3p-2, 2p+1}, (-3p-2)( 2p+1))$. 
Let $D'$ be a disk bounded by $c'$. 
Then $T_{-3p-2, 2p+1} \cup D'$ can be isotoped so that $T_{-3p-2, 2p+1}$ lies on $\Sigma$, 
$D' \cap \Sigma$ consists of a single arc, 
which intersects $T_{-3p-2, 2p+1}$ in $4p+3$ points in the same direction.  
Now taking the mirror image of $T_{-3p-2, 2p+1} \cup D'$, 
we obtain $T_{3p+2, 2p+1} \cup D$ with $\partial D = c$; 
$D \cap \Sigma$ consists of a single arc, 
and $D$ intersects $T_{3p+2, 2p+1}$ in $4p+3$ points in the same direction. 
Then $c$ is a seiferter for $(T_{3p+2, 2p+1}, (3p+2)( 2p+1))$. 
Since $(4p+3)^2 > 2(3p+2)(2p+1)$, 
Theorem~\ref{Tpq} shows that $K(3p+2, 2p+1; 4p+3, n)$ is an L-space knot for all integers $n$. 
Let us show that $K(3p+2, 2p+1; 4p+3, n)$ is hyperbolic if $|n| >3$. 
Figure~47 in \cite{DMM2}, 
together with \cite[Lemma~5.6(1)]{DMM1}, 
shows that $c'$ becomes a degenerate fiber in $T_{-3p-2, 2p+1}((-3p-2)( 2p+1))$, 
and so $c$ becomes a degenerate fiber in $T_{3p+2, 2p+1}((3p+2)( 2p+1))$. 
Apply the same argument as above to obtain the desired result.

\begin{figure}[h]
\begin{center}
\includegraphics[width=0.6\linewidth]{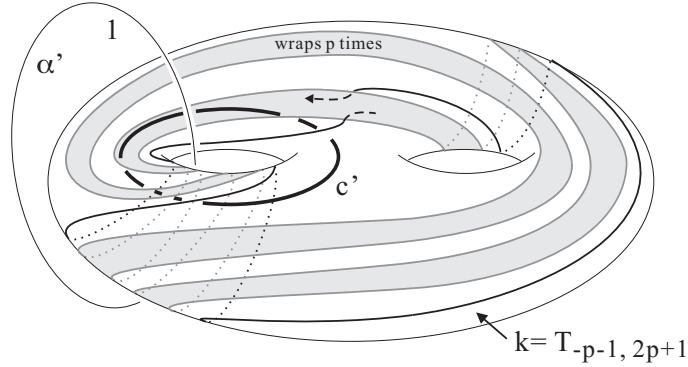}
\caption{A surgery description of $T_{-3p-2, 2p+1}$ and a seiferter $c'$}
\label{sporcfiberA}
\end{center}
\end{figure}

\medskip

$\bullet$ $K(2p+3, 2p+1; 2p+2, n)$ $(p > 0)$. \ 
We follow the argument in \cite[Section~6]{DMM2}; 
as above we need to take the mirror image at the end. 
Take a torus knot $k = T_{-3p-2, 3}$ on a genus two Heegaard surface of $S^3$, 
unknotted circles $\alpha'$ and  $c'$ as in Figure~\ref{BergeVI}(i). 
Then $(-2)$--twist along $\alpha'$ converts the torus knot $T_{-3p-2, 3}$ into a Berge knot $\mathrm{VI}[p]$. 
Lemma~6.1 in \cite{DMM2} shows that $c'$, 
the image of $c'$ after the $(-2)$--twist along $\alpha'$, 
is a seiferter for a lens space surgery $(\mathrm{VI}[p], -8p^2 -16p-7)$; 
the linking number between $c'$ and $\mathrm{VI}[p]$ is $2p+2$. 
Now we show that $1$--twist along $c'$ (after $(-2)$--twist along $\alpha'$) converts   
$(\mathrm{VI}[p], -8p^2 -16p-7)$ into $(T_{-2p-1, 2p+3}, (-2p-1)(2p+3))$. 
Note that $c'$ remains a seiferter for $(T_{-2p-1, 2p+3}, (-2p-1)(2p+3))$. 
Since the linking number between $c'$ and $\mathrm{VI}[p]$ is $2p+2$, 
the surgery slope $-8p^2 -16p-7$ becomes $-8p^2 -16p-7 + (2p+2)^2 = (-2p-1)(2p+3)$. 
Let us observe that the knot obtained from $\mathrm{VI}[p]$ by $1$--twist along $c'$, 
which has a surgery description given by Figure~\ref{BergeVI}(i), 
is $T_{-2p-1, 2p+3}$. 
The surgeries described in Figure~\ref{BergeVI}(i) can be realized by the following two successive twistings: 
$1$--twist along an annulus cobounded by $c'$ and $\alpha'$ (cf.\cite[Definition 2.32]{DMM1}), 
and $(-1)$--twist along $\alpha'$. 
The annulus twist converts $k = T_{-3p-2, 3}$ into $k' = T_{-2p-1, 2}$ as shown in Figure~\ref{BergeVI}(ii). 
Then $(-1)$--twist along $\alpha'$ changes $k' = T_{-2p-1, 2}$ into $T_{-2p-1, 2p+3}$, 
which lies on the genus one Heegaard surface $\Sigma$. 
Let $D'$ be a disk bounded by $c'$. 
Then $D'$ can be slightly isotoped so that 
$D' \cap \Sigma$ consists of a single arc, 
which intersects $T_{-2p-1, 2p+3}$ in $2p+2$ points in the same direction; 
see Figure~\ref{BergeVI}(ii). 
Now taking the mirror image of $T_{-2p-1, 2p+3} \cup D'$, 
we obtain $T_{2p+1, 2p+3} \cup D$ with $\partial D = c$; 
$D \cap \Sigma$ consists of a single arc, 
which intersects $T_{2p+1, 2p+3}$ in $2p+2$ points in the same direction. 
Then $c$ is a seiferter for $(T_{2p+1, 2p+3}, (2p+1)(2p+3)) = (T_{2p+3, 2p+1}, (2p+3)(2p+1))$. 
Theorem~\ref{Tpq} shows that $K(2p+3, 2p+1; 2p+2, n)$ is an L-space knot for any integer $n \ge -1$. 
The hyperbolicity of knots $K(2p+3, 2p+1; 2p+2, n)$ for $|n| > 3$ follows from the same argument as above, 
in which we refer to Figure~33 instead of Figure~47 in \cite{DMM2}. 

\begin{figure}[h]
\begin{center}
\includegraphics[width=1.0\linewidth]{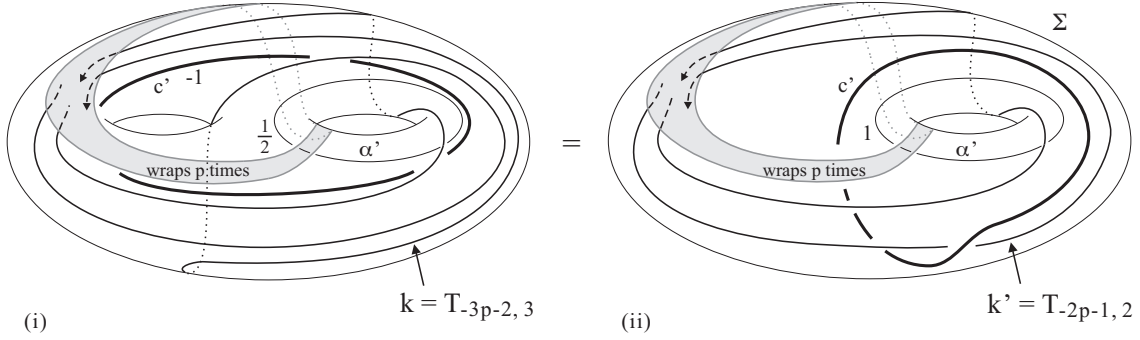}
\caption{Surgery descriptions of $T_{-2p-1, 2p+3}$ and a seiferter $c'$}
\label{BergeVI}
\end{center}
\end{figure}

\hspace*{\fill} $\square$(Theorem~\ref{twisted torus knot})

\medskip

\noindent
\textit{Proof of Corollary~\ref{L-space twisted torus knot}.}
Given any torus knot $T_{p, q}$ $(p > q \ge 2)$, 
let us take an unknotted circle $c = c^+_{p. q}$ in $S^3 - T_{p, q}$ (Figure~\ref{cpq_plus}(i)). 
Then as shown in the proof of Theorem~\ref{twisted torus knot}, 
$n$--twist along $c$ converts $T_{p, q}$ into the twisted torus knot 
$K(p, q; p+q, n)$, which is an L-space knot for all integers $n$ and hyperbolic if $|n| > 3$. 

The last assertion of Corollary~\ref{L-space twisted torus knot} follows from 
Claim~\ref{distinct twisted torus knots} below. 
Thus the unknotted circle $c$ satisfies the required property in Corollary~\ref{L-space twisted torus knot}. 
\hspace*{\fill} $\square$(Corollary~\ref{L-space twisted torus knot})

\medskip

\begin{claim}
\label{distinct twisted torus knots}
$\{ K(p, q; p+q, n) \}_{|n| > 3}$ is a set of mutually distinct hyperbolic knots. 
\end{claim}

\noindent
\textit{Proof of Claim~\ref{distinct twisted torus knots}.}
Recall that $c^+_{p, q}$ is a seiferter for $(T_{p, q}, pq)$ and the linking number between $ c^+_{p, q}$ and 
$T_{p, q}$ is $p+q$. 
Thus $n$--twist along $c^+_{p, q}$ changes $(T_{p, q}, pq)$ to a Seifert surgery $(K(p, q; p + q, n), pq + n(p+q)^2)$. 
Note that $K(p, q; p + q, n)(pq + n(p+q)^2)$ is a Seifert fiber space over $S^2$ with at most three exceptional fibers 
of indices $p$, $q$ and $|n|$, see the proof of Theorem~\ref{twisted torus knot}. 

Assume that $K(p, q; p + q, n)$ is isotopic to $K(p, q; p+q, n')$ for some integers 
$n, n'$ with $|n|, |n'| > 3$. 
Then $pq + n(p+q)^2$--, and $pq + n'(p+q)^2$--surgeries on the hyperbolic knot $K(p, q; p + q, n)$ yield 
Seifert fiber spaces. 
Hence $|pq + n(p+q)^2 - (pq + n'(p+q)^2)| = |(n-n')(p+q)^2| \le 8$ by \cite[Theorem~1.2]{LacMey}. 
Since $p+q \ge 5$, we have $n = n'$. 
This completes a proof. 
(In the above argument, we can apply \cite[Theorem~8.1]{Agol} which gives the bound $10$ instead of $8$.)
\hspace*{\fill} $\square$(Claim~\ref{distinct twisted torus knots})

\bigskip

\section{L-space twisted Berge knots}
\label{twisted Berge knots}

In this section we prove Theorem~\ref{Berge}  
using Theorem~\ref{Tpq} and observations in \cite{DMM1, DMM2}. 

Berge \cite{Berge2} gave twelve infinite families of knots which admit lens space surgeries. 
These knots are referred to as \textit{Berge knots} of types $(\mathrm{I})$--$(\mathrm{XII})$ and conjectured to comprise all knots with lens space surgeries. 
Recall that a Berge knot of type $(\mathrm{I})$ is a torus knot and that of $(\mathrm{II})$ is a cable of a torus knot, 
henceforth we consider Berge knots of types $(\mathrm{III})$--$(\mathrm{XII})$. 

\medskip

$\bullet$ 
Berge knots of types $(\mathrm{III})$--$(\mathrm{VI})$.

Suppose that $K$ is a Berge knot of type $(\mathrm{III})$, $(\mathrm{IV})$, $(\mathrm{V})$ or $(\mathrm{VI})$. 
Then we have an unknotted solid torus $V$ containing $K$ in its interior such that $V(K; m)$ is a solid torus \cite{Berge2, DMM2}, 
and hence the core $c$ of the solid torus $W = S^3 - \mathrm{int}V$ is a seiferter for $(K, m)$ and 
$(K_n, m_n)$ is also a lens space. 
If $K_n(m_n)$ is not an L-space, 
then it is $S^2 \times S^1$ and $(K_n, m_n) = (O, 0)$ (\cite[Theorem~8.1]{GabaiIII}). 
Now let us exclude this possibility. 
First we note that $V(K_n, m_n) \cong V(K; m)$ for all integers $n$ and 
$H_1(V(K_n; m_n)) \cong \mathbb{Z} \oplus \mathbb{Z}_{(m_n, \omega)}$ \cite[Lemma~3.3]{Go1}, 
where $\omega$ is the winding number of $K$ in $V$, 
i.e. the linking number between $K_n$ and $c$. 
Since $V(K_n; m_n) \cong S^1 \times D^2$, 
$K_n$ is a $0$ or $1$--bridge braid in $V$ \cite{Gabai_solidtorus}, 
hence $\omega \ge 2$. 
This then implies that $m_n \ne 0$. 
Hence $(K_n, m_n)$ is an L-space knot for all integers $n$. 

\medskip

$\bullet$ 
Berge knots of types $(\mathrm{VII})$, $(\mathrm{VIII})$. 

Let $g_1$ and $g_2$ be simple closed curves embedded in a genus two Heegaard surface 
$F$ of $S^3$ and $c$ an unknot in $S^3$ as in Figure~\ref{Berge78}. 

\begin{figure}[h]
\begin{center}
\includegraphics[width=0.45\linewidth]{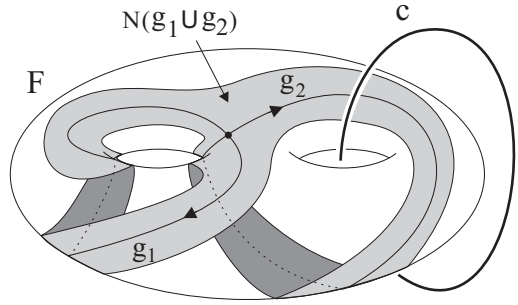}
\caption{A regular neighborhood $N(g_1 \cup g_2)$ of $g_1 \cup g_2$ in $F$ and an unknotted circle $c$}
\label{Berge78}
\end{center}
\end{figure}

Take a regular neighborhood $N(g_1 \cup g_2)$ of $g_1 \cup g_2$ in $F$, 
which is a once punctured torus. 
Then the curve $\partial N(g_1 \cup g_2)$ becomes a trefoil knot after $(-1)$--twist along $c$, 
and the figure-eight knot after $1$--twist along $c$. 
Let $k$ be a knot in $N(g_1 \cup g_2)$ representing $a[g_1] + b[g_2] \in H_1(N(g_1 \cup g_2))$, 
where $a$ and $b$ are coprime integers. 
The we see that $k$ is a torus knot $T_{a+b, -a}$. 
The Berge knot $K$ of type $(\mathrm{VII})$ (resp. $(\mathrm{VIII})$) is obtained from $T_{a+b, -a}$ by $(-1)$--twist 
(resp. $1$--twist) along $c$. 
As shown in \cite[Lemma~4.6]{DMM1}, 
$T_{a+b, -a} \cup c$ is isotopic to $T_{a+b, -a} \cup c^+_{a+b, -a}$ and 
Berge knot of type $(\mathrm{VII})$ is $K(a+b, -a; |b|, -1)$, 
that of type $(\mathrm{VIII})$ is $K(a+b, -a; |b|, 1)$; 
see the proof of Theorem~\ref{twisted torus knot}. 
(Here we extend the notation $K(p, q; r, n)$ for twisted torus knots 
in an obvious fashion to include the case where $p, q$ are possibly negative integers.)
We assume $|a|, |b| \ge 2$, 
for otherwise $K(a+b, -a, |b|, \pm 1)$ is a torus knot. 
Furthermore, if $|a+b| = 1$, 
then $T_{a+b, -a} \cup c = T_{\pm 1, -a} \cup c$ is a torus link $T_{2, 2b}$ or $T_{2, -2b}$, 
and  $K(a+b, -a; |b|, \pm 1)$ is a torus knot,  
so we assume $|a+b| > 1$.  
Let $K_n$ be a knot obtained from the Berge knot $K$ by $n$--twist along $c$, 
i.e. $K_n = K(a+b, -a; |b|, n+\varepsilon)$; $\varepsilon = -1$ if $K$ is of type $(\mathrm{VII})$, 
$\varepsilon = 1$ if $K$ is of type $(\mathrm{VIII})$. 
If $a(a+b) < 0$ (i.e. $-a(a+b) > 0$), 
then by Theorem~\ref{twisted torus knot} $K_n$ is an L-space knot for any integer $n$. 
If $a(a+b) > 0$ (i.e. $-a(a+b) < 0$), 
Theorem~\ref{twisted torus knot} shows that the mirror image $K(a+b, a; |b|, -n-\varepsilon)$ of $K_n$ is an L-space knot if $-n-\varepsilon \ge -1$, i.e. $n \le 1 - \varepsilon$. 
Hence $K_n$ is an L-space knot for any integer $n \le 1 - \varepsilon$.  

\medskip 

$\bullet$ Berge knots of types $(\mathrm{IX})$--$(\mathrm{XII})$. 

These knots are often called \textit{sporadic} knots and we denote them by 
$Spor{\bf a}[p]$, $Spor{\bf b}[p]$, $Spor{\bf c}[p]$ and $Spor{\bf d}[p]$ $(p \ge 0)$, respectively.  
It is easy to see that 
$Spor{\bf a}[0]$ and $Spor{\bf b}[0]$ are trivial knots, 
$Spor{\bf c}[0] =T_{-3, 4}$ and $Spor{\bf c}[0] = T_{-5, 3}$. 
Thus we may assume $p > 0$ for $Spor{\bf \chi}[p]$ $({\bf \chi} = {\bf a}, {\bf b}, {\bf c}, {\bf d})$. 
Furthermore, we observe that 
$Spor{\bf a}[1]$ is obtained from $T_{3, 2}$ by $1$--twist along the seiferter $c = c^+_{3, 2}$; 
see Figure~\ref{T32_5}. 
Hence following Example~\ref{trefoil_c6} a knot $K_n$ obtained from $Spor{\bf a}[1]$ by 
$n$--twist along $c$ is an L-space knot for any integer $n$. 
Thus we may assume $p > 1$ for $Spor{\bf a}[p]$. 

As shown in Example~\ref{cable} the lens space surgery $(Spor{\bf a}[p], 22p^2+9p+1)$ is obtained from 
$(C_{6p+1, p}(T_{3, 2}), p(6p+1))$ by $1$--twist along the seiferter $c$, 
and $n$--twist along $c$ converts $C_{6p+1, p}(T_{3, 2})$ into an L-space knot for all integers $n$.  
Hence $n$--twist changes $Spor{\bf a}[p]$ to an L-space knots for all integers $n$. 

The proof of Theorem~\ref{twisted torus knot} shows that 
the lens space surgery $(Spor{\bf b}[p], 22p^2 + 13p+2)$ is obtained from $(T_{3p+1, 2p+1}, (3p+1)(2p+1))$ 
by $1$--twist along the seiferter $c$, hence $K_n$ obtained from $Spor{\bf b}[p]$ by $n$--twist along $c$ 
is $K(3p+1, 2p+1; 4p+1, n+1)$. 
By Theorem~\ref{twisted torus knot} $K_n$ is an L-space knot for all integers $n$. 
Similarly  $(Spor{\bf c}[p], -22p^2-31p-11)$ is obtained from  $(T_{-3p-2, 2p+1}, (-3p-2)( 2p+1))$ 
by $(-1)$--twist along $c'$, 
and $K_n$ obtained from $Spor{\bf c}[p]$ by $n$--twist along $c'$ 
is $K(-3p-2, 2p+1; 4p+4, n-1)$, 
which is the mirror image of $K(3p+2, 2p+1; 4p+4, -n+1)$. 
Theorem~\ref{twisted torus knot} shows that $K(3p+2, 2p+1; 4p+4, -n+1)$ is an L-space knot for any integer $n$, 
and thus $K_n$ is an L-space knot for all integers $n$. 

Finally, let us consider a Berge knot $Spor{\bf d}[p]$ $(p \ge 0)$. 
Proposition~8.8 in \cite{DMM2} shows that 
the lens space surgery $(Spor{\bf d}[p], -22p^2 -35p-14)$ has a seiferter $c'$ 
such that the linking number between $c'$ and $Spor{\bf d}[p]$ is $4p+3$ 
and $1$--twist along $c'$ converts $(Spor{\bf d}[p], -22p^2 -35p-14)$ into 
$(C_{-6p-5, p+1}(T_{-3, 2}), (-6p-5)(p+1))$ for which $c'$ is a seiferter. 
Let $K_n$ be a knot obtained from $Spor{\bf d}[p]$ by $n$--twist along $c'$, 
i.e. obtained from $C_{-6p-5, p+1}(T_{-3, 2})$ by $(n-1)$--twist along $c'$. 
Now we take the mirror image of $C_{-6p-5, p+1}(T_{-3, 2}) \cup c'$ to obtain 
a link $C_{6p+5, p+1}(T_{3, 2}) \cup c$. 
Then $c$ is a seiferter for $(C_{6p+5, p+1}(T_{3, 2}), (6p+5)(p+1))$ and 
$K_n$ is the mirror image of the knot obtained from $C_{6p+5, p+1}(T_{3, 2})$ by 
$(-n)$--twist along $c$. 
Since $(4p+3)^2 \ge 2(6p+5)(p+1)$, Theorem~\ref{Tpq} shows that 
$K_n$ is an L-space knot for all integers $n$.  

\medskip

Let us show that $K_n$ is a hyperbolic knot except for at most four integers $n$. 
Following \cite[Theorem 5.10]{DMM1} it is sufficient to observe that $K \cup c$ is a hyperbolic link. 
Suppose that $K$ is a Berge knot of type $(\mathrm{III})$, $(\mathrm{IV})$, $(\mathrm{V})$ or $(\mathrm{VI})$. 
Then as mentioned above, 
$V(K; m)$ is a solid torus, where $V = S^3 - \mathrm{int}N(c)$. 
By \cite[Theorem~3.2]{Berge1} $V- \mathrm{int}N(K)$ is atoroidal. 
If $V- \mathrm{int}N(K)$ is not hyperbolic, 
then it is Seifert fibered and $K$ is a torus knot; see \cite[Lemma~3.3]{DMM1}. 
This contradicts the assumption. 
Hence $K \cup c$ is a hyperbolic link. 

If $K$ is of type $(\mathrm{VII})$ or $(\mathrm{VIII})$, 
then $K \cup c \cong T_{a+b, -a} \cup c_{a+b, -a}^+$ is a hyperbolic link; see the proof of Theorem~\ref{twisted torus knot}. 

Assume that $K$ is of type $(\mathrm{IX})$, i.e. $K = Spor{\bf a}[p]$. 
Then as shown in the proof of Example~\ref{cable} $K \cup c$ is a hyperbolic link. 
In the case where $K$ is of type $(\mathrm{X})$ or $(\mathrm{XI})$, 
i.e. $K = Spor{\bf b}[p]$ or $Spor{\bf c}[p]$, 
then it follows from the proof of Theorem~\ref{twisted torus knot} that $K \cup c$ is a hyperbolic link. 
The argument in the proof of Example~\ref{cable} shows $K \cup c$ is a hyperbolic link for 
type $(\mathrm{XII})$ Berge knot $K = Spor{\bf d}[p]$; 
we refer to Figure~53 instead of Figure~41. 

This completes a proof of Theorem~\ref{Berge}. 
\hspace*{\fill} $\square$(Theorem~\ref{Berge})

\bigskip

\section{L-space twisted unknots}
\label{twisted unknots}

In \cite{DMM1} we introduced ``$m$--move" to find seiferters for a given Seifert surgery. 
In particular, 
$m$--move is effectively used in \cite[Theorem~6.21]{DMM1} to show 
that $(O, m)$ has infinitely many seiferters for each integer $m$. 
Among them there are infinitely many seiferters $c$ such that 
$(m, 0)$--surgery on $O \cup c$ is an L-space; see Remark~\ref{non Lspace seiferter}. 

Let us take a trivial knot $c_{m, p}$ in $S^3 - O$ as illustrated in Figure~\ref{Om_seiferter}, 
where $p$ is an odd integer with $|p| \ge 3$.  

\begin{figure}[htbp]
\begin{center}
\includegraphics[width=0.6\linewidth]{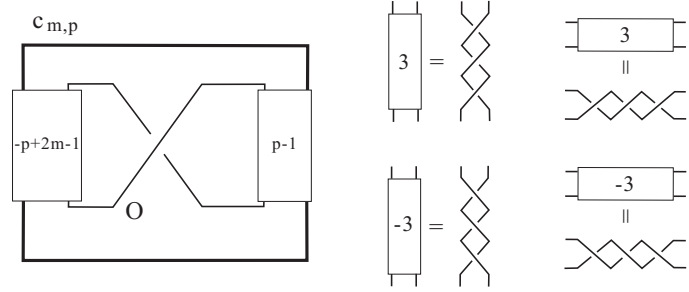}
\caption{$O \cup c_{m, p}$; a vertical (resp.\ horizontal) box with integer $n$ denotes a vertical (resp.\ horizontal) stack of $n$ crossings.}
\label{Om_seiferter}
\end{center}
\end{figure}

Then as shown in \cite[Theorem~6.21]{DMM1}
$c_{m, p}$ is a seiferter for $(O, m)$ such that $O \cup c_{m, p}$ is a hyperbolic link in $S^3$ if $p \ne 2m \pm 1$.  
Denote by $K_{m, p, n}$ and $m_{p, n}$ the images of $O$ and $m$ after $n$--twist along $c_{m, p}$. 
Now we investigate $K_{m, p, n}(m_{p, n})$ using branched coverings and Montesinos trick \cite{Montesinos, Mon}. 
Figure~\ref{Om_seiferter1} (b) shows that $K_{m, p, n}(m_{p, n})$ has an involution with axis L for any integer $n$. 
Taking the quotient by this involution, 
we obtain a $2$--fold branched cover $\pi : K_{m, p, n}(m_{p, n}) \to S^3$ branched along $L'$ which is the quotient of L; 
see Figure~\ref{Om_seiferter1}(c).  
As shown in Figure~\ref{Om_seiferter1}(d) 
$L'$ can be isotoped to a Montesinos link $M(-n/(mn+1),\ (-p+1)/2p,\ (-p+2m +1)/ (-2p+4m))$. 
Hence by \cite{Montesinos} $K_{m, p, n}(m_{p, n})$, 
which is the $2$--fold branched cover branched along the Montesinos link $L'$,  
is a Seifert fiber space 
$$S^2(\frac{-n}{mn+1},\ \frac{-p+1}{2p},\ \frac{p-2m -1}{2p-4m}).$$ 
The image $\pi(c_{m, p})$ is an arc $\tau$ whose ends lie in $L'$; see Figure~\ref{Om_seiferter1}(c) and (d). 
It follows from \cite[Lemma~3.2]{DEMM} that $c_{m, p}$ is a seiferter for $(K_{m, p, n}, m_{p, n})$; 
in case of $n = 0$, $c_{m, p}$ is a seiferter for $(O, m)$. 
In the following, 
the image of $c_{m, p}$ after $n$--twist along itself 
is denoted by the same symbol. 

\begin{figure}[htbp]
\begin{center}
\includegraphics[width=1.0\linewidth]{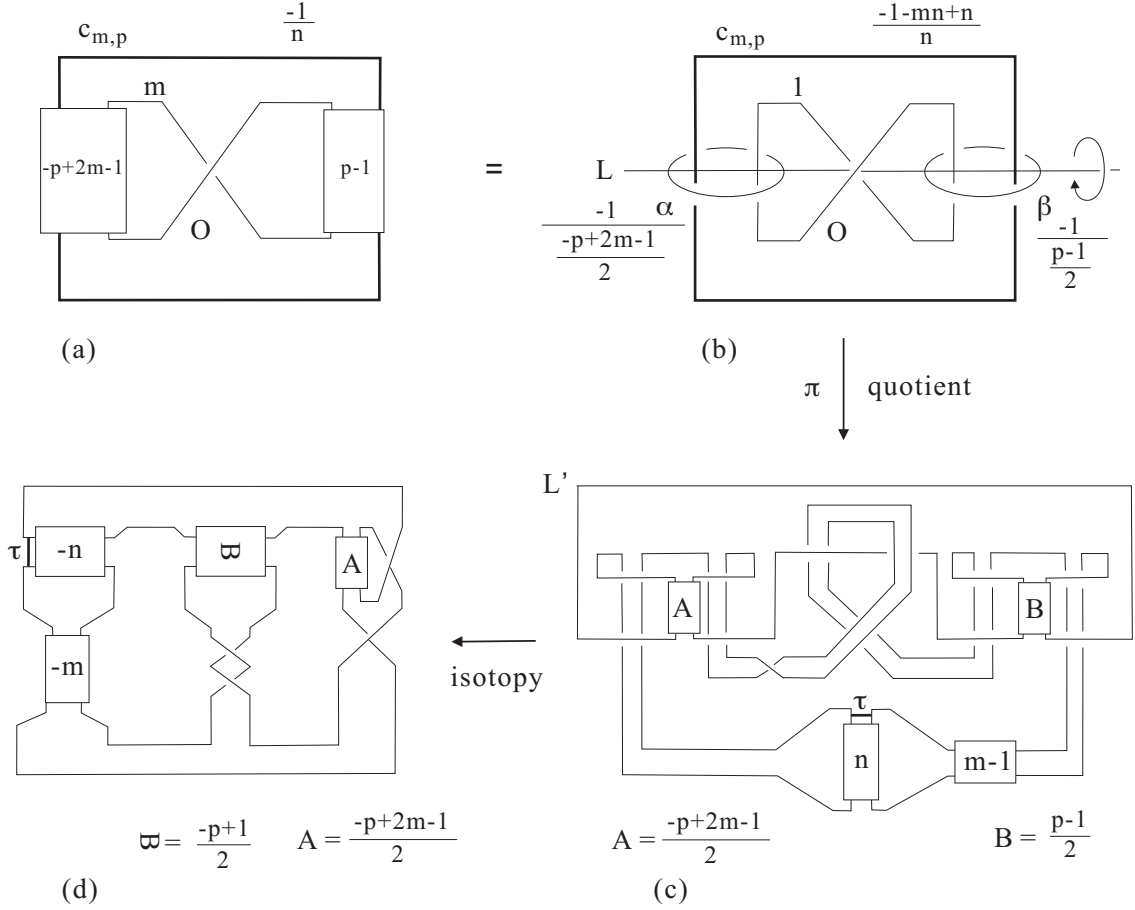}
\caption{$K_{m, p, n}(m_{p, n})$ is the two-fold branched cover of $S^3$ branched along $L'$.}
\label{Om_seiferter1}
\end{center}
\end{figure}

In what follows assume $m \le 0$ and $p \ge 3$. 

\begin{proposition} 
\label{Kn(mn)L-space}
Assume that $m \le 0,\ p \ge 3$. 
\begin{enumerate}
\item
$(K_{m, p, n}, m_{p, n})$ is an L-space surgery except when $(m, n) = (0, 0)$. 
If $(m, n) = (0, 0)$, 
then $(K_{m, p, n}, m_{p, n}) = (O, 0)$ and $K_{m, p, n}(m_{p, n}) = O(0) \cong S^2 \times S^1$. 

\item
$K_{m, p, n}$ is a nontrivial knot if $n \ne 0$. 

\item
$\{ K_{m, p, n} \}_{|n| > 1}$ is a set of mutually distinct hyperbolic L-space knots. 

\end{enumerate}

\end{proposition}

\noindent
\textit{Proof of Proposition~\ref{Kn(mn)L-space}.}
We note here that the linking number between $c_{m, p}$ and $O$ is $p-m$. 

$(1)$ Assume first that $m = 0$. 
Then $K_{m, p, n}(m_{p, n})$ is a lens space 
$S^2(-n,\  (-p+1)/2p,\ (p-1)/2p) = S^2(-n-1,\ (p+1)/2p,\ (p-1)/2p)$, 
which is $S^2 \times S^1$ if and only if $n = 0$ (Claim~\ref{S2xS1}).  
Hence $K_{m, p, n}(m_{p, n})$ is an L-space except when $n = 0$. 

Next assume $m = -1$. 
Then $K_{m, p, n}(m_{p, n}) 
= S^2(-n/(-n+1),\ (-p+1)/2p,\ (p+1)/ (2p+4)) = 
S^2(n/(n-1),\ (-p+1)/2p,\ (p+1)/ (2p+4))$. 
If $n = 0$ or $2$, 
$K_{m, p, n}(m_{p, n})$  is a lens space, but it is not  $S^2 \times S^1$, 
because $m_{p, n} = -1 + n(m-p)^2 = -1 + n(p+1)^2 \ne 0$.  
If $n = 1$, $K_{m, p, n}(m_{p, n})$ is a connected sum of two lens spaces,  
and thus an L-space. 
Suppose that $n \ne 0, 1, 2$. 
In case of $n < 0$, 
we have $0 < n/(n-1) < 1$ and 
$K_{m, p, n}(m_{p, n}) 
= S^2(n/(n-1),\ (-p+1)/2p,\ (p+1)/ (2p+4)) 
= S^2(-1,\ n/(n-1),\ (p+1)/2p,\ (p+1)/ (2p+4))$. 
Note that $(p+1)/2p + (p+1)/ (2p+4) 
= 1/2 + 1/2p + 1/2 - 1/ (2p+4)  = 1 + 1/2p - 1/(2p+4)$. 
Since $p \ge 3$, 
we have $2p+4 > 2p > 0$, 
and hence $1/2p - 1/(2p+4) > 0$. 
It follows that $(p+1)/2p +(p+1)/ (2p+4) = 1 + 1/2p - 1/(2p+4) > 1$. 
Then Lemma~\ref{any r}(2) shows that $K_{m, p, n}(m_{p, n})$ is an L-space. 
If $n > 2$, then $1 < n/(n-1) < 2$ and 
$K_{m, p, n}(m_{p, n}) = S^2(n/(n-1),\ (-p+1)/2p,\ (p+1)/ (2p+4)) 
= S^2(1/(n-1),\ (p+1)/2p,\ (p+1)/ (2p+4))$. 
Since $0 < 1/(n-1),\ (p+1)/2p,\ (p-2m -1)/ (2p-4m) < 1$, 
$K_{m, p, n}(m_{p, n})$ is an L-space 
by Theorem~\ref{SeifertLspace}(1). 

Assume that $m = -2$. 
Then $K_{m, p, n}(m_{p, n}) 
= S^2(-n/(-2n+1),\ (-p+1)/2p,\ (p+3)/ (2p+8)) 
= S^2(n/(2n-1),\ (-p+1)/2p,\ (p+3)/ (2p+8))$. 
If $n = 0, 1$, 
then $K_{m, p, n}(m_{p, n})$ is a lens space, but it is not $S^2 \times S^1$, 
because $m_{p, n} = -2 + n(m-p)^2 = -2 + n(p+2)^2 \ne 0$.  
Otherwise, 
$0 < n/(2n-1) < 1$ and 
$K_{m, p, n}(m_{p, n}) 
= S^2(n/(2n-1),\ (-p+1)/2p,\ (p+3)/ (2p+8)) = 
S^2(-1,\ n/(2n-1),\ (p+1)/2p,\ (p+3)/ (2p+8))$. 
Since $(p+1)/2p + (p+3)/ (2p+8) 
= 1/2 + 1/2p +1/2 -1/(2p+8) = 1 + 1/2p - 1/(2p+8) > 1$, 
$K_{m, p, n}(m_{p, n})$ is an L-space by Lemma~\ref{any r}(2). 

Finally assume that $m \le -3$. 
then $K_{m, p, n}(m_{p, n}) 
= S^2(-n/(mn+1),\ (-p+1)/2p,\ (p-2m -1)/ (2p-4m)) 
= S^2(-1,\ -n/(mn+1),\ (p+1)/2p,\ (p-2m -1)/ (2p-4m))$. 
If $n= 0$, then $K_{m, p, n}(m_{p, n})$ is a lens space, but it is not $S^2 \times S^1$, 
because $m_{p, n} = m + n(m-p)^2 = m \le -3$. 
Assume $n \ne 0$. 
Then by the assumption $p \ge 3, m \le -3$ we have $0 < -n/(mn+1) < 1$, 
$0 < (p+1)/2p < 1$ and 
$0 < (p-2m -1)/ (2p-4m) = 1/2 - 1/(2p-4m) < 1$. 
Since $(p+1)/2p + (p-2m -1)/ (2p-4m) 
= 1/2 + 1/2p + 1/2 - 1/(2p-4m) =1 + 1/2p - 1/(2p-4m) > 1$, 
Lemma~\ref{any r}(2) shows that $K_{m, p, n}(m_{p, n})$ is an L-space.  

\medskip

$(2)$ 
Since $m \le 0$ and $p \ge 3$, 
$p \ne 2m \pm 1$, 
hence $O \cup c_{m, p}$ is a hyperbolic link; 
see \cite[Theorem 6.21]{DMM1}. 
Then $K_{m, p, n}$ is nontrivial for any $n \ne 0$ \cite{KMS, Ma}. 

\medskip

$(3)$ 
By $(1)$ $K_{m, p, n}$ is an L-space knot. 
Since $O \cup c_{m, p}$ is a hyperbolic link, 
the hyperbolicity of $K_{m, p, n}$ for $|n| > 1$ follows from \cite{AMM2, GLtoroidal, MM3}.  
Thus $K_{m, p, n}$ $(|n| > 1)$ is a hyperbolic L-space knot. 
Let us choose $c_{m, p}$ and then apply $n$--twist along $c_{m, p}$ to obtain a knot $K_{m, p, n}$. 
It remains to show that $K_{m, p, n}$ and $K_{m, p, n'}$ are distinct knots. 
Suppose that $K_{m, p, n}$ and $K_{m, p, n'}$ are isotopic for some integers $n$ and $n'$ with $|n|, |n'| > 1$.  
Then $(m+ n(p-m)^2)$--, and $(m+n'(p-m)^2)$--surgeries on $K_{m, p, n} = K_{m, p, n'}$ produce small Seifert fiber spaces, 
where $p-m \ge 3$. 
(Note that since $|n| > 1$, $mn+1$ cannot be zero.)
Since $K_{m, p, n}$ is a hyperbolic knot, 
Lackenby and Meyerhoff \cite[Theorem~ 1.2] {LacMey} prove that the distance $|m+n(p-m)^2 - (m+n'(p-m)^2|$ between above two non-hyperbolic surgeries 
is at most $8$. 
Hence $|(n-n')(p-m)^2| \le 8$, 
which implies $n = n'$ because $p-m \ge 3$. 
\hspace*{\fill} $\square$(Proposition~\ref{Kn(mn)L-space})

\medskip

Next we investigate link types of $O \cup c_{m, p}$.  

\begin{proposition}
\label{infinity}
Let $c_{m,p}$ and $c_{m', p'}$ be seiferters for $(O,  m)$ and $(O, m')$, respectively. 
Suppose that $m, m' \le 0,\ p, p' \ge 3$. 
\begin{enumerate}
\item 
If $p-m \ne p' - m'$, 
then $O \cup c_{m, p}$ and $O \cup c_{m', p'}$ are not isotopic. 
In particular, 
if $p \ne p'$, 
then $O \cup c_{m, p}$ and $O \cup c_{m, p'}$ are not isotopic. 

\item
If $p-m = p' - m'$, 
then $O \cup c_{m, p}$ and $O \cup c_{m', p'}$ are not isotopic provided that 
$|m - m'| > 3$. 
\end{enumerate}

\end{proposition}

\noindent
\textit{Proof of Proposition~\ref{infinity}.}\ 
$(1)$ Note that the linking number between $c_{m, p}$ and $O$ is $p-m$. 
Hence if $O \cup c_{m, p}$ is isotopic to $O \cup c_{m', p'}$ as ordered links, 
then we have $p-m = p' - m'$. 

\medskip

$(2)$ Since $p \ne 2m \pm 1$ and $p' \ne 2m' \pm 1$, 
both $O \cup c_{m, p}$ and $O \cup c_{m', p'}$ are hyperbolic links \cite{DMM1}. 
Recall that $c_{m, p}$ is a seiferter for $(O, m)$ and $c_{m', p'}$ is a seiferter for $(O, m')$. 
Suppose that $O \cup c_{m, p}$ and $O \cup c_{m', p'}$ are isotopic. 
Then $c_{m, p}$ is a seiferter for $(O, m')$ as well. 
Let $V$ be the solid torus $S^3 - \mathrm{int}N(c_{m, p})$, which contains $O$ in its interior. 
Note that $m$--surgery of $V$ along $O$ yields a Seifert fiber space over the disk with two exceptional fibers 
of indices $2p, 2p-4m$, 
and $m'$--surgery of $V$ along $O$ yields a Seifert fiber space over the disk with two exceptional fibers 
of indices $2p', 2p'-4m'$. 
Since these Seifert fiber spaces contain essential annuli,  
Gordon and Wu \cite[Corollary~1.2]{GW} show that $|m - m'| \le 3$. 
\hspace*{\fill} $\square$(Proposition~\ref{infinity})

\medskip

Theorem~\ref{L-space twisted unknot} follows from Propositions~\ref{Kn(mn)L-space} and \ref{infinity}. 
\hspace*{\fill} $\square$(Theorem~\ref{L-space twisted unknot})

\begin{remark}
\label{non Lspace seiferter}
For each seiferter $c_{m, p}$ $(m \le 0,\ p \ge 3)$, 
$M_{c_{m, p}}(O, m)$ is an L-space. 
In fact, 
$M_{c_{m, p}}(O, m)$, 
which is the limit of $K_{m, p, n}(m_{p, n})$ when $|n|$ tends to $\infty$ $($Remark~\ref{limit}$)$, 
is $S^2(-1/m, (-p+1)/2p, (p-2m-1)/(2p-4m)) = S^2(-1, -1/m, (p+1)/2p, (p-2m-1)/(2p-4m))$. 
If $m = -1, 0$, then $M_{c_{m, p}}(O, m)$ is an L-space $($Claim~\ref{rc}$)$. 
If $m < -1$, since $(p+1)/2p + (p-2m-1)/(2p-4m) = 1 + 1/2p - 1/(2p-4m) > 1$, 
$M_{c_{m, p}}(O, m)$ is an L-space. 

On the other hand, for instance, 
$M_{c_{3, 3}}(O, 3)$ is not an L-space. 
Indeed, $M_{c_{3, 3}}(O, 3) = S^2(-1/3, -1/3, 2/3) = S^2(-2, 2/3, 2/3, 2/3)$, 
and taking $k = 2,\ a = 1$ in Theorem~\ref{SeifertLspace}(3), 
we have $(1-2/3, 1-2/3, 1-2/3) = (1/3, 1/3 , 1/3) < (1/2, 1/2, 1/2)$. 
Thus $M_{c_{3, 3}}(O, 3)$ is not an L-space. 
\end{remark}

\bigskip

\section{Hyperbolic, L-space knots with tunnel number greater than one}
\label{tunnel}

The purpose in this section is to exhibit infinitely many hyperbolic L-space knots with tunnel number greater than one (Theorem~\ref{tunnel2}).  
In \cite{EJMM} Eudave-Mu\~noz, Jasso and Miyazaki and the author gave 
Seifert fibered surgeries which do not arise from primitive/Seifert-fibered construction \cite{D}. 

Let us take unknotted circles $c_a$ and $c_b$ in $S^3 - T_{3, 2}$ as illustrated by Figure~\ref{trefoil_annularpair}. 
Then as shown in \cite{EJMM} $\{ c_a, c_b \}$ is a \textit{pair of seiferters} for $(T_{3, 2}, 7)$, 
i.e. $c_a$ and $c_b$ become fibers simultaneously in some Seifert fibration of 
$T_{3, 2}(7)$. 

\begin{figure}[h]
\begin{center}
\includegraphics[width=0.3\linewidth]{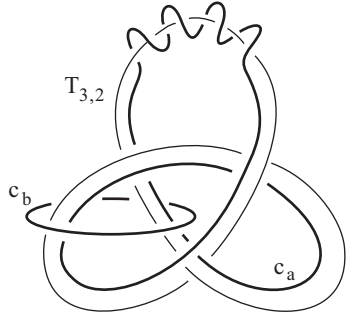}
\caption{$\{c_a, c_b\}$ is a pair of seiferters for $(T_{3,2}, 7)$.}
\label{trefoil_annularpair}
\end{center}
\end{figure}

Note that the pair $\{c_a,\ c_b \}$ forms the $(4, 2)$--torus link in $S^3$. 
Hence $(-1)$--twist along $c_a$ converts $c_a \cup c_b$ into the $(-4, 2)$--torus link. 
Then we can successively apply $1$--twist along $c_b$ to obtain $(4, 2)$--torus link $c_a \cup c_b$. 
We denote the images of $c_a,\ c_b$ under twistings along these components by the same symbols $c_a,\ c_b$,  respectively.

Let $K_{n, 0}$ be a knot obtained from $T_{3, 2}$ after the sequence of twistings:  
$$(c_a, (-1)\mathrm{-twist}) \to  (c_b, 1\mathrm{-twist}) \to (c_a, n\mathrm{-twist})$$
Then $K_{n, 0} = K(2, -n, 1, 0)$ in \cite[Proposition~4.11]{EJMM}. 
See Figure~\ref{T32_tunnel2}. 

Similarly, let $K_{0, n}$ be a knot obtained from $T_{3, 2}$ after the sequence of twistings: 
$$(c_a, (-1)\mathrm{-twist}) \to  (c_b, n+1\mathrm{-twist})$$
Then $K_{0, n} = K(2, 0, 1, -n)$ in \cite[Proposition~4.11]{EJMM}.  
See Figure~\ref{T32_tunnel2}. 

\begin{figure}[h]
\begin{center}
\includegraphics[width=0.5\linewidth]{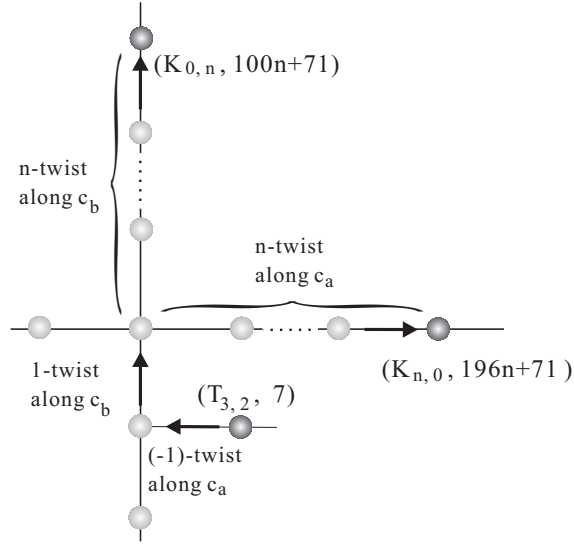}
\caption{Seifert surgeries $(K_{n, 0}, 196n+71)$ and $(K_{0, n}, 100n+71)$; 
each vertex corresponds to a Seifert surgery and each edge corresponds to a single twist along a seiferter.}
\label{T32_tunnel2}
\end{center}
\end{figure}

Theorem~\ref{tunnel2} follows from Theorem~\ref{Kn0K0n} below. 

\begin{theorem}
\label{Kn0K0n}
\begin{enumerate}
\item
$\{ K_{n, 0} \}_{n \in \mathbb{Z}}$ is a set of mutually distinct hyperbolic L-space knots with tunnel number two. 

\item
$\{ K_{0, n} \}_{n \in \mathbb{Z} \setminus \{ -1 \}}$ is a set of mutually distinct hyperbolic L-space knots with tunnel number two. 
\end{enumerate}
\end{theorem}

\noindent
\textit{Proof of Theorem~\ref{Kn0K0n}.}
We begin by recalling the following result which is a combination of Propositions~3.2,\ 3.7 and 3.11 in \cite{EJMM}. 

\begin{lemma}
\label{Kn0K0n_tunnel}
\begin{enumerate}
\item $K_{n, 0}$ is a hyperbolic knot with tunnel number two, 
and $K_{n, 0}(196n+71)$ is a Seifert fiber space 
$S^2((11n+4)/(14n+5),\ -2/7,\ 1/2)$. 

\item
$K_{0, n}$ is a hyperbolic knot with tunnel number two if $n \ne -1$, 
and $K_{0, n}(100n+71)$ is a Seifert fiber space 
$S^2(-(3n+2)/(10n+7),\ 4/5,\ 1/2)$. 

\end{enumerate}
\end{lemma}

\begin{lemma}
\label{Kn0K0n_distinct}
\begin{enumerate}
\item
If $K_{n, 0}$ and $K_{n', 0}$ are isotopic, then $n=n'$.
\item
If $K_{0, n}$ and $K_{0, n'}$ are isotopic, then $n = n'$. 
\end{enumerate}
\end{lemma}

\noindent
\textit{Proof of Lemma~\ref{Kn0K0n_distinct}.}
$(1)$ Suppose that $K_{n, 0}$ is isotopic to $K_{n', 0}$. 
Then $K_{n, 0}(196n+71)$ and $K_{n, 0}(196n'+71)$ are both Seifert fiber spaces. 
Since $K_{n, 0}$ is hyperbolic,  
Theorem~1.2 in \cite{LacMey} implies that $|196n+71-(196n'+71)| = |196(n-n')| \le 8$. 
Hence we have $n = n'$. 
$(2)$ follows in a similar fashion. 
\hspace*{\fill} $\square$(Lemma~\ref{Kn0K0n_distinct})

\medskip

Let us prove that $K_{n, 0}$ and $K_{0, n}$ are L-space knots for any integer $n$. 

\begin{lemma}
\label{tunnel2L-space}
\begin{enumerate}
\item
$K_{n, 0}(196n+71)$ is an L-space for any integer $n$. 
\item
$K_{0, n}(100n+71)$ is an L-space for any integer $n$. 
\end{enumerate}
\end{lemma}

\noindent
\textit{Proof of Lemma~\ref{tunnel2L-space}.}
$(1)$ Note that $K_{n, 0}(196n+71) = S^2((11n+4)/(14n+5),\ -2/7,\ 1/2) = S^2(-1, (11n+4)/(14n+5),\ 5/7,\ 1/2)$. 
Since $0 < (11n+4)/(14n+5) < 1$ for any $n \in \mathbb{Z}$ and $5/7 + 1/2 \ge 1$, 
Lemma~\ref{any r}(2) shows that $K_{n, 0}(196n+71)$ is an L-space for any integer $n$. 
This proves (1). 

$(2)$ As above first we note that 
$K_{0, n}(100n+71) = S^2(-(3n+2)/(10n+7),\ 4/5,\ 1/2) = S^2(-1,\ (7n+5) /(10n+7),\ 4/5,\ 1/2)$. 
Since $0 < (7n+5) /(10n+7) < 1$ for any $n \in \mathbb{Z}$ and $4/5 + 1/2 \ge 1$, 
Lemma~\ref{any r}(2) shows that $K_{0, n}(100n+71)$ is an L-space for any integer $n$. 
\hspace*{\fill} $\square$(Lemma~\ref{tunnel2L-space})

\medskip

Now Theorem~\ref{Kn0K0n} follows from Lemmas~\ref{Kn0K0n_tunnel}, \ref{Kn0K0n_distinct} and \ref{tunnel2L-space}. 
\hspace*{\fill} $\square$(Theorem~\ref{Kn0K0n})

\begin{question}
\label{large tunnel}
Does there exist a hyperbolic L-space knot with tunnel number greater than two? 
More generally, for a given integer $p$, 
does there exist a hyperbolic L-space knot with tunnel number greater than $p$?
\end{question}

\bigskip

\section{Questions}
\label{questions}

\subsection{Characterization of twistings which yield infinitely many L-space knots.}

For knots $K$ with Seifert surgery $(K, m)$, 
Theorems~\ref{tree small}, \ref{tree projective}, \ref{tree degenerate} and Corollary~\ref{Tpq} 
characterize seiferters which enjoy the desired property in Question~\ref{motivation}. 

The next proposition, 
which is essentially shown in \cite{Hedden, HLV}, 
describes yet another example of twistings which yield infinitely many L-space knots. 

\begin{proposition}[\textbf{L-space twisted satellite knots}]
\label{twisted satellite}
Let $k$ be a nontrivial knot with L-space surgery $(k, 2g-1)$, where $g$ denotes the genus of $k$, 
and $K$ a satellite knot of $k$ which lies in 
$V = N(k)$ with winding number $w$. 
Suppose that $V(K; m)$ is a solid torus for some integer $m \ge w^2(2g-1)$. 
Let $c$ be the boundary of a meridian disk of $V$, 
and $K_n$ a knot obtained from $K$ by $n$--twist along $c$. 
Then $K_n$ is an L-space knot for any $n \ge 0$.  See Figure~\ref{satellite}. 
\end{proposition}

\noindent
\textit{Proof of Proposition~\ref{twisted satellite}.} 
Recall that $K_n(m+ nw^2) = k((m+ nw^2)/w^2) = k(m/w^2 + n)$ \cite{Go1}. 
Since $k(2g-1)$ is an L-space and $m/w^2 \ge 2g-1$, 
\cite[Proposition 9.6]{OS4} ensures that $k(m/w^2 + n)$ is also an L-space if $n \ge 0$.  
Hence $K_n$ is an L-space knot provided $n \ge 0$. 
\hspace*{\fill} $\square$(Proposition~\ref{twisted satellite})

\medskip

\begin{figure}[h]
\begin{center}
\includegraphics[width=0.3\linewidth]{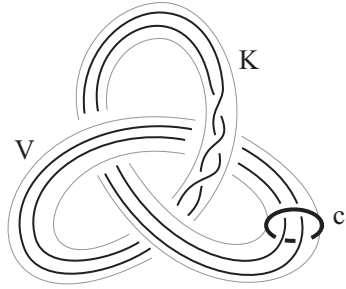}
\caption{$K_n$ is a knot obtained from $K$ by $n$--twist along $c$.}
\label{satellite}
\end{center}
\end{figure}

\begin{remark}
\label{satellite_twist}
\begin{enumerate}
\item
In Proposition~\ref{twisted satellite}, 
the knot $K$ in the solid torus $V$ is required to have a cosmetic surgery: $V(K; m) \cong S^1 \times D^2$. 
The cosmetic surgery of the solid torus is well-understood by \cite{Gabai_solidtorus, Berge1}. 
\item
Twisting operation described in Proposition~\ref{twisted satellite} can be applied only for satellite knots and 
resulting knots after the twistings are also satellite knots.  
\item
In Proposition~\ref{twisted satellite}, the knot $k$ is assumed to be nontrivial. 
If $k$ is a trivial knot in $S^3$, 
then $K(m) = (S^3 - \mathrm{int}V) \cup V(K; m)$ is a lens space, hence $(K, m)$ is an L-space surgery. 
It is easy to see that $c$ is a seiferter for $(K, m)$. 
\end{enumerate}
\end{remark}

For further possibility, 
weaken a condition of seiferter to obtain a notion of ``pseudo-seiferter" as follows. 

\begin{definition}[\textbf{pseudo-seiferter}]
Let $(K, m)$ be a Seifert surgery.
A knot $c$ in $S^3 -N(K)$ is called a \textit{pseudo-seiferter} for $(K, m)$
if $c$ satisfies (1) and (2) below.
\begin{enumerate}
\item $c$ is a trivial knot in $S^3$.
\item $c$ becomes a ``cable" of a fiber in
a Seifert fibration of $K(m)$ and the preferred longitude $\lambda$ of $c$ in $S^3$ becomes the cabling slope of $c$ in $K(m)$. 
\end{enumerate}
\end{definition}

We do not know if a pseudo-seiferter exists, 
but if $(K, m)$ admits a pseudo-seiferter, it behaves like a seiferter in the following sense. 
Let $V$ be a fibered tubular neighborhood of a fiber $t$ and $c$ is a cable in $V$. 
Then the result of a surgery (corresponding to $n$--twist) on $c$ of $V$ is again a solid torus, 
and this surgery is reduced to a surgery on the fiber $t$ which is a core of $V$. 
Hence $K_n(m_n)$ is a (possibly degenerate) Seifert fiber space. 
This suggests that a pseudo-seiferter is also a candidate for an unknotted circle described in Question~\ref{motivation}. 

We would like to ask the following question for non-satellite knots. 

\begin{question}
\label{twist=seiferter}
Let $K$ be a non-satellite knot and $K_n$ a knot obtained from $K$ by $n$--twist along an unknotted circle $c$ in $S^3 - K$. 
Suppose that the twist family $\{ K_n \}$ contains infinitely many L-space knots. 

\begin{enumerate}
\item
Does $K$ admit a Seifert surgery $(K, m)$ for which $c$ is a seiferter? 

\item
Does $K$ admit a Seifert surgery $(K, m)$ for which $c$ is a seiferter or a pseudo-seiferter? 
\end{enumerate}
\end{question}

\medskip

\subsection{L-space knots and strong invertibility.}

A knot is said to be \textit{strongly invertible} 
if there exists an orientation preserving involution of $S^3$ which fixes the knot setwise and reverses orientation.
Known L-space knots are strongly invertible, 
so it is natural to ask: 

\begin{problem}[Watson]
\label{invertible}
Are L-space knots strongly invertible?
\end{problem}

In \cite{DMM1} an ``asymmetric seiferter" defined below is essentially used to find Seifert fibered sugary on 
knots with no symmetry.  

\begin{definition}[\textbf{asymmetric seiferter}]
A seiferter $c$ for a Seifert surgery $(K, m)$ 
is said to be \textit{symmetric} 
if we have an orientation preserving diffeomorphism 
$f : S^3 \to S^3$ of finite order with $f(K) = K, f(c) = c$;
otherwise, $c$ is called an \textit{asymmetric seiferter}. 
\end{definition}

Combining \cite[Theorem~7.3]{DMM1} and Theorem~\ref{tree small}, we obtain: 

\begin{proposition}
\label{asymmetric}
Let $(K, m)$ be a Seifert fibered surgery on a non-satellite knot with an asymmetric seiferter $c$ 
which becomes an exceptional fiber. 
Suppose that $M_c(K, m)$ is an L-space. 
Then there is a constant $N$ such that $K_n$,
a knot obtained from $K$ by $n$--twist along $c$, is a hyperbolic L-space knot with no symmetry 
for any $n \le N$ or $n \ge N$. 
\end{proposition}

If $c$ is a seiferter for $(T_{p, q}, pq)$ which becomes a degenerate fiber in $T_{p, q}(pq)$, 
then $c$ is a meridian of $T_{p, q}$ or $T_{p, q} \cup c$ is a hyperbolic link in $S^3$; 
see \cite[Theorem~3.19(3)]{DMM1}. 
Hence the argument in the proof of Theorem~7.3 in \cite{DMM1} and Theorem~\ref{tree degenerate}(2) 
enable us to show: 

\begin{proposition}
\label{asymmetric2}
If $c$ is an asymmetric seiferter for $(T_{p, q}, pq)$ which becomes a degenerate fiber in $T_{p, q}(pq)$, 
then there is a constant $N$ such that $K_n$ is a hyperbolic L-space knot with no symmetry 
for any $n \le N$ or $n \ge N$. 
\end{proposition}

\medskip

For the asymmetric seiferter $c = c'_1$ for $(K, m) = (P(-3,3,5), 1)$ given in \cite[Lemma~7.5]{DMM1}, 
$M_c(K, m)$ is not an L-space and $c$ does not satisfy the hypothesis of Proposition~\ref{asymmetric}. 

\begin{question}
\label{asymmetric seiferter}
Does there exist an asymmetric seiferter described in Propositions~\ref{asymmetric} and \ref{asymmetric2}? 
\end{question}


\bigskip

\begin{thebibliography}{99}


\bibitem{Agol}
I. Agol; 
Bounds on exceptional Dehn filling, 
Geom.\ Topol.\ \textbf{4} (2000), 431--449. 

\bibitem{AMM2}
M. A\"it Nouh, D. Matignon and K. Motegi; 
Geometric types of twisted knots, 
Ann.\ Math.\ Blaise Pascal \textbf{13} (2006), 31--85. 

\bibitem{BM}
K. Baker and A. Moore; 
Montesinos knots, Hopf plumbings, and L-space surgeries, 
preprint 2014. 

\bibitem{BePe}
R. Benedetti and C. Petronio; 
Lectures on hyperbolic geometry, 
Universitext, Springer-Verlag, 1992. 

\bibitem{Berge1} 
J. Berge; 
The knots in $D^2 \times S^1$ which have nontrivial Dehn surgeries 
that yield $D^2 \times S^1$, 
Topology Appl.\ \textbf{38} (1991), 1--19.

\bibitem{Berge2} 
J. Berge; 
Some knots with surgeries yielding lens spaces, 
unpublished manuscript.

\bibitem{BoileauPorti}
M. Boileau and J. Porti; 
Geometrization of $3$-orbifolds of cyclic type,  
Ast\'erisque \textbf{272} (2001), 208pp. 

\bibitem{BGW}
S. Boyer, C. McA. Gordon and L. Watson; 
On L-spaces and left-orderable fundamental groups, 
Math.\ Ann.\ \textbf{356} (2013), 1213--1245. 

\bibitem{BRW} 
S. Boyer, D. Rolfsen and B. Wiest; 
Orderable 3-manifold groups, 
Ann.\ Inst.\ Fourier \textbf{55} (2005), 243--288. 

\bibitem{D} 
J. Dean; 
Small Seifert-fibered Dehn surgery on hyperbolic knots, 
Algebr.\ Geom.\ Topol.\ \textbf{3} (2003), 435--472.  


\bibitem{DMM2}
A. Deruelle, K. Miyazaki and K. Motegi; 
Networking Seifert Surgeries on Knots II:  Berge's lens surgeries, 
Topology Appl.\ \textbf{156} (2009), 1083--1113. 

\bibitem{DEMM}
A. Deruelle, M. Eudave-Mu\~noz, K. Miyazaki and K. Motegi; 
Networking Seifert surgeries on knots IV: seiferters and branched coverings, 
Contemp.\ Math.\ Amer.\ Math.\ Soc.\ \textbf{597} (2013), 235--262. 

\bibitem{DMM1} 
A. Deruelle, K. Miyazaki and K. Motegi; 
Networking Seifert surgeries on knots, 
Mem.\ Amer.\ Math.\ Soc.\ \textbf{217} (2012), no. 1021, viii+130. 

\bibitem{EHN}
D. Eisenbud, U. Hirsch and W. Neumann; 
Transverse foliations on Seifert bundles and self-homeomorphisms of the circle, 
Comment.\ Math.\ Helv.\ \textbf{56} (1981), 638--660. 

\bibitem{EM2} 
M. Eudave-Mu\~noz; 
On hyperbolic knots with Seifert fibered Dehn surgeries, 
Topology Appl.\ \textbf{121} (2002), 119--141. 

\bibitem{EJMM}
M. Eudave-Mu\~noz, E. Jasso, K. Miyazaki and K. Motegi; 
Seifert fibered surgeries on strongly invertible knots without primitive/Seifert positions, 
to appear in Topology Appl..  
DOI:10.1016/j.topol.2015.05.016

\bibitem{GabaiIII}
D. Gabai; 
Foliations and the topology of $3$--manifolds. III, 
J.\ Diff.\ Geom.\ \textbf{26} (1987), 479--536

\bibitem{Gabai_solidtorus}
D. Gabai; 
Surgery on knots in solid tori, 
Topology \textbf{28} (1989), 1--6. 

\bibitem{Ghig}
P. Ghiggini; 
Knot Floer homology detects genus-one fibred knots, 
Amer.\ J.\ Math.\ \textbf{130} (2008), 1151--1169. 

\bibitem{Go1} 
C.McA. Gordon;  
Dehn surgery and satellite knots, 
Trans.\ Amer.\ Math.\ Soc.\ \textbf{275} (1983), 687--708.  

\bibitem{GLtoroidal}
C.McA. Gordon and J. Luecke; 
Non-integral toroidal Dehn surgeries, 
Comm.\ Anal.\ Geom.\ \textbf{12} (2004), 417--485. 

\bibitem{GW}
C.McA. Gordon and Y.-Q. Wu; 
Annular Dehn fillings, 
Comment.\ Math.\ Helv.\ \textbf{75} (2000), 430--456. 

\bibitem{Greene}
J. E. Greene; 
L-space surgeries, genus bound, and the cabling conjecture, 
J.\ Diff.\ Geom.\ \textbf{100} (2015), 491--506. 

\bibitem{Hat}
A. E. Hatcher; 
Notes on basic $3$-manifold topology, 
freely available at \texttt{http://www.math.cornell.edu/hatcher}. 

\bibitem{Hedden}
M. Hedden; 
On knot Floer homology and cabling II,  
Int.\ Math.\ Res.\ Not.\ IMRN, (12):2248--2274,
2009.

\bibitem{HLV}
J. Hom, T. Lidman and F. Vafaee; 
Berge-Gabai knots and L-space satellite operations, 
Algebr.\ Geom.\ Topol.\ \textbf{14} (2014) 3745--3763.  


\bibitem{Ja}
W. Jaco;
Lectures on three manifold topology, 
CBMS Regional Conference Series in Math., vol. 43,
Amer.\ Math.\ Soc., 1980.


\bibitem{JN}
M. Jankins and W. Neumann; 
Rotation number and products of circle homomorphisms, 
Math.\ Ann.\ \textbf{271} (1985), 381--400. 

\bibitem{Juh}
A. Juh\'asz; 
Floer homology and surface decompositions, 
Geom.\ Topol.\ \textbf{12} (2008), 299--350. 

\bibitem{KMS}
M. Kouno, K. Motegi and T. Shibuya;   
Twisting and knot types, 
J.\ Math.\ Soc.\ Japan \textbf{44} (1992), 199-216. 

\bibitem{Kr}
D. Krcatovich; 
The reduced Knot Floer complex, 
preprint 2014. 

\bibitem{LacMey}
M. Lackenby and R. Meyerhoff, 
The maximal number of exceptional Dehn surgeries, 
Invent.\ Math.\  \textbf{191} (2013), 341--382. 

\bibitem{LM}
T. Lidman and A. Moore; 
Pretzel knots with L-space surgeries, 
to appear in Mich.\ Math.\  J. 

\bibitem{LS}
P. Lisca and A. Stipsicz;
Ozsv\'ath-Szab\'o invariants and tight contact 3-manifolds, III.,  
J.\ Symplectic Geom.\ \textbf{5}, (2007), 357?384. 

\bibitem{Ma} Y. Mathieu; 
Unknotting, knotting by twists on disks and property (P) for 
knots in $S^3$, Knots 90 (ed. by Kawauchi), 
Proc. 1990 Osaka Conf. on Knot Theory and Related Topics, 
de Gruyter (1992), 93-102

\bibitem{MM1}
K. Miyazaki and K. Motegi; 
Seifert fibred manifolds and Dehn surgery, 
Topology \textbf{36} (1997), 579--603. 

\bibitem{MM3}
K. Miyazaki and K. Motegi; 
Seifert fibered manifolds and Dehn surgery III, 
Comm.\ Anal.\ Geom.\ \textbf{7} (1999) 551--582. 

\bibitem{Montesinos} 
J. M. Montesinos; 
Variedades de Seifert que son recubridores c\'iclicos ramificados de dos hojas, 
Bol.\ Soc.\ Mat.\ Mexicana \textbf{18} (1973), 1--32. 

\bibitem{Mon}
J. M. Montesinos; 
Surgery on links and double branched coverings of $S^3$, 
Knots, groups, and $3$--manifolds (Papers dedicated to the memory of R.H.Fox),
Ann. Math. Studies, 84, Princeton Univ. Press 1975, 227--260,  

\bibitem{MoSa}
K. Morimoto and M. Sakuma; 
On unknotting tunnels for knots, 
Math.\ Ann.\ \textbf{289} (1991), 143--167. 

\bibitem{MT}
K. Motegi and K. Tohki; 
On L-space twisted torus knots, 
in preparation. 

\bibitem{Na}
R. Naimi; 
Foliations transverse to fibers of Seifert manifolds, 
Comment.\ Math.\ Helv.\ \textbf{69} (1994), 155--162. 

\bibitem{Ni}
Y. Ni; 
Knot Floer homology detects fibred knots, 
Invent.\ Math.\  \textbf{170} (2007) 577--608. 

\bibitem{Ni2}
Y. Ni; 
Erratum: Knot Floer homology detects fibred knots, 
Invent.\ Math.\ \textbf{177} (2009), 235--238. 

\bibitem{OS2}
P. Ozsv\'ath and Z. Szab\'o; 
Holomorphic disks and three-manifold invariants:\ properties and applications, 
Ann.\ of Math.\ 
\textbf{159} (2004), 1159--1245. 

\bibitem{OSgenus}
P. Ozsv\'ath and Z. Szab\'o; 
Holomorphic disks and genus bounds, 
Geom.\ Topol.\ \textbf{8} (2004), 311--334. 

\bibitem{OS3} 
P. Ozsv\'ath and Z. Szab\'o; 
On knot Floer homology and lens space surgeries, 
Topology \textbf{44} (2005), 1281--1300. 

\bibitem{OS4} 
P. Ozsv\'ath and Z. Szab\'o; 
Knot Floer homology and rational surgeries, 
Algebr.\ Geom.\ Topol.\ \textbf{11} (2011), 1--68. 

\bibitem{PetPorti}
C. Petronio and J. Porti;
Negatively oriented ideal triangulations and a proof of Thurston's hyperbolic Dehn filling theorem, 
Expo.\ Math.\ \textbf{18} (2000), 1--35. 

\bibitem{Szabo}
Z. Szab\'o; 
Lecture Notes on Heegaard Floer Homology, 
LAS/Park City Mathematics Series, Volume 15 (2006), 199--228. 

\bibitem{T1}
W. P. Thurston; 
The geometry and topology of $3$-manifolds, 
Lecture notes, Princeton University, 1979. 

\bibitem{T2}
W. P. Thurston; 
Three dimensional manifolds, Kleinian groups and hyperbolic geometry, 
Bull.\ Amer.\ Math.\ Soc.\ \textbf{6} (1982), 357--381. 

\bibitem{V}
F. Vafaee; 
On the knot Floer homology of twisted torus knots, 
Int.\ Math.\ Res.\ Notices.\ \textbf{2015-15} (2015), 6516--6537. 


\end{thebibliography}
\end{document}